\documentclass[11pt]{article}

\usepackage{graphicx}
\usepackage{amsfonts}
\usepackage{amsmath}
\usepackage{amsthm}
\usepackage[mathcal]{euscript}      % Define tensors
\usepackage{bm}                     % Define vectors
\usepackage{color}
\usepackage{epstopdf}
\usepackage{amssymb}
\usepackage{mathrsfs,mathdots}
\usepackage{lscape}

\def \ep{\hbox{ }\hfill$\Box$}

\usepackage[paper=a4paper,dvips,top=3cm,left=2.4cm,right=2.4cm,
    foot=1.4cm,bottom=3cm]{geometry}

\begin{document}

\title{Positive Definiteness of Paired Symmetric Tensors and Elasticity Tensors}% and Strong Ellipticity}
\author{Zheng-Hai Huang\footnote{%
    School of Mathematics, Tianjin University, Tianjin 300354, P.R. China ({\tt huangzhenghai@tju.edu.cn}).
    This author was supported by the National Natural Science Foundation of China (Grant No. 11431002).}
  \and Liqun Qi\footnote{%
    Department of Applied Mathematics, The Hong Kong Polytechnic University,
    Hung Hom, Kowloon, Hong Kong ({\tt maqilq@polyu.edu.hk}).
    This author's work was partially supported by the Hong Kong Research Grant Council
    (Grant No. PolyU  501913, 15302114, 15300715 and 15301716).}
    }

\date{\today}
\maketitle

\begin{abstract}
In this paper, we consider higher order paired symmetric tensors and strongly paired symmetric tensors.  Elasticity tensors and higher order elasticity tensors in solid mechanics are strongly paired symmetric tensors. A (strongly) paired symmetric tensor is said to be positive definite if the homogeneous polynomial defined by it is positive definite.   Positive definiteness of elasticity and higher order elasticity tensors is strong ellipticity in solid mechanics, which plays an important role in nonlinear elasticity theory.   We mainly investigate positive definiteness of fourth order three dimensional and sixth order three dimensional (strongly) paired symmetric tensors. We first show that the concerned (strongly) paired symmetric tensor is positive definite if and only if its smallest $M$-eigenvalue is positive. Second, we propose several necessary and sufficient conditions under which the concerned (strongly) paired symmetric tensor is positive definite. Third, we study the conditions under which the homogeneous polynomial defined by a fourth order three dimensional or sixth order three dimensional (strongly) paired symmetric tensor can be written as a sum of squares of polynomials, and further, propose several necessary and/or sufficient conditions to judge whether the concerned (strongly) paired symmetric tensors are positive definite or not. Fourth, by using semidefinite relaxation we propose a sequential semidefinite programming method to compute the smallest $M$-eigenvalue of a fourth order three dimensional (strongly) paired symmetric tensor, by which we can check positive definiteness of the concerned tensor. The preliminary numerical results confirm our theoretical findings.

\end{abstract}

\newtheorem{Theorem}{Theorem}[section]
\newtheorem{Definition}[Theorem]{Definition}
\newtheorem{Lemma}[Theorem]{Lemma}
\newtheorem{Corollary}[Theorem]{Corollary}
\newtheorem{Proposition}[Theorem]{Proposition}
\newtheorem{Example}[Theorem]{Example}%[section]
\newtheorem{Remark}[Theorem]{Remark}%[section]

% LaTeX definitions
\renewcommand{\hat}[1]{\widehat{#1}}
\renewcommand{\tilde}[1]{\widetilde{#1}}
\renewcommand{\bar}[1]{\overline{#1}}
\newcommand{\REAL}{\mathbb{R}}
\newcommand{\COMPLEX}{\mathbb{C}}
\newcommand{\SPHERE}{\mathbb{S}^2}
\newcommand{\diff}{\,\mathrm{d}}
\newcommand{\st}{\mathrm{s.t.}}
\newcommand{\T}{\top}
\newcommand{\vt}[1]{{\bf #1}}%{\bm{#1}}
\newcommand{\x}{\vt{x}}
\newcommand{\y}{\vt{y}}
\newcommand{\z}{\vt{z}}
\newcommand{\Ten}{\mathcal{T}}
\newcommand{\A}{\mathcal{A}}
\newcommand{\RESULTANT}{\mathrm{Res}}

\newpage
\section{Introduction}

For any positive integers $m$ and $n$, an $m$-th order $n$ dimensional real tensor can be denoted by
$$
\mathscr{A}:=(a_{i_1i_2\cdots i_m}),\; \mbox{\rm where}\;a_{i_1i_2\cdots i_m}\in \mathbb{R}\; \mbox{\rm for all}\; i_j\in [n]\; \mbox{\rm with}\; j\in [m],
$$
here $[l]:=\{1,2,\ldots,l\}$ for any positive integer $l$. We use $\mathbb{T}_{m,n}$ to denote the set of all $m$-th order $n$ dimensional real tensors. For any $\mathscr{A}=(a_{i_1i_2\cdots i_m}), \mathscr{B}=(b_{i_1i_2\cdots i_m})\in \mathbb{T}_{m,n}$, we will use the inner product defined by
$$
\langle\mathscr{A},\mathscr{B}\rangle=\sum_{i_1,\cdots,i_m=1}^na_{i_1i_2\cdots i_m}b_{i_1i_2\cdots i_m}
$$
and the Hilbert-Schmidt norm defined by
$$
\|\mathscr{A}\|_{HS}=\sqrt{\langle\mathscr{A},\mathscr{A}\rangle}=\left(\sum_{i_1,\cdots,i_m=1}^na_{i_1i_2\cdots i_m}^2\right)^{1/2}.
$$
A tensor $\mathscr{A}\in \mathbb{T}_{m,n}$ is said to be symmetric if its entries are invariant under any permutation of indices $\{i_1,i_2,\ldots,i_m\}$. For any $x=(x_1,\ldots,x_n)\in \mathbb{R}^n$, a tensor $\mathscr{A}=(a_{i_1i_2\cdots i_m})\in \mathbb{T}_{m,n}$ defines a homogeneous polynomial by
\begin{eqnarray}\label{e-func-general}
p(x)=\sum_{i_1,\cdots,i_m=1}^na_{i_1i_2\cdots i_m}x_{i_1}x_{i_2}\cdots x_{i_m}.
\end{eqnarray}
We denote the degree of a polynomial $p$ by deg$(p)$. The polynomial $p$ defined by (\ref{e-func-general}) is said to be positive semidefinite if $p(x)\geq 0$ holds for all $x\in \mathbb{R}^n$; and $p$ is said to be positive definite if $p(x)>0$ holds for all $x\in \mathbb{R}^n\backslash \{0\}$. Obviously, for nonzero tensors, $m$ being an even integer is a necessity for positive
semidefiniteness. A tensor is said to positive (semidefinite) definite if its corresponding homogeneous polynomial is positive (semidefinite) definite. Positive (semidefiniteness) definiteness of the polynomial is very important in many areas, which is related to the Hilbert seventeenth problem \cite{Reznick-95}.

For an arbitrary $2m$-th order $n$ dimensional tensor denoted by $\mathscr{A}=(a_{i_1j_1i_2j_2\cdots i_mj_m})$, its indices can be divided into $m$ adjacent blocks $\{i_1j_1\}, \ldots, \{i_mj_m\}$. If entries of $\mathscr{A}$ are invariant under any permutation of indices in every block $\{i_lj_l\}$ for $l\in [m]$, i.e.,
$$
a_{i_1j_1i_2j_2\cdots i_mj_m}=a_{j_1i_1i_2j_2\cdots i_mj_m}=a_{i_1j_1j_2i_2\cdots i_mj_m}=\cdots =a_{i_1j_1i_2j_2\cdots j_mi_m},
$$
then $\mathscr{A}$ is called a {\it $2m$-th order $n$ dimensional paired symmetric tensor}. It is well known that the most important representatives for fourth order three dimensional paired symmetric tensors are: the piezooptical tensor, the second order electrooptical effect (Kerr effect), electrostriction and second order magnetostriction \cite{Ha-07}. Furthermore, if a paired symmetric tensor additionally satisfies block symmetry, i.e.,
$$
a_{i_1j_1i_2j_2\cdots i_mj_m}=a_{i_2j_2i_3j_3\cdots i_{1}j_{1}}=a_{i_3j_3i_4j_4\cdots i_{2}j_{2}}=\cdots =a_{i_mj_mi_{1}j_{1}\cdots i_{m-1}j_{m-1}},
$$
then $\mathscr{A}$ is called a {\it $2m$-th order $n$ dimensional strongly paired symmetric tensor}. The most significant representative of fourth order three dimensional strongly paired symmetric tensors is the elasticity tensor, in which the pairwise permutability is based on the reversibility of mechanical deformation
work \cite{Ha-07}. In the elasticity tensor and higher order elasticity tensor \cite{Hiki-81,JCAJ-15,TB-64}, every entry $a_{i_1j_1i_2j_2\cdots i_mj_m}$ is called an $m$-th order elastic constant, which is an important quantity in studies of elasticity theory. Positive definiteness of the elasticity tensor is called  {\it strong ellipticity}, which plays an important role in elasticity theory and has been studied extensively (see, for example, \cite{CG-10,HDQ-09,KS-75,KS-77,Mu-13,QDH-09,Ro-90,Sf-11,SS-83,St-13,WA-96,ZR-16}).

In this paper, we consider positive definiteness of higher order three dimensional (strongly) paired symmetric tensors. For simplicity of symbols, we only investigate some properties of fourth order three dimensional and sixth order three dimensional (strongly) paired symmetric tensors, mainly in positive definiteness of the concerned tensors. The results we obtained can be similarly extended to the case of more higher order (strongly) paired symmetric tensors. Some basic properties of fourth order three dimensional and sixth order three dimensional (strongly) paired symmetric tensors are given in the next section.

Eigenvalues of higher order tensors, introduced by Qi \cite{Qi-05} and Lim \cite{Lim-05}, have been studied extensively in the recent years \cite{CS-13,HHLQ-13,Qi-07,QL-17}. The concept of $M$-eigenvalue for fourth order paired symmetric tensor was introduced in \cite{HDQ-09,QDH-09} and further studied in \cite{WQZ-09}. In Section 3, we extend the concept of $M$-eigenvalue to sixth order three dimensional (strongly) paired symmetric tensors and bi-block symmetric tensors, and further discuss some related properties. In particular, we show that a sixth order three dimensional (strongly) paired symmetric tensor is positive definite if and only if its smallest $M$-eigenvalue is positive.

Positive definiteness of polynomials has been being an important issue in many areas, which has been discussed extensively. In \cite{HDQ-09,NZ-16,QDH-09}, the authors studied positive definiteness conditions of fourth order paired symmetric tensors, which plays an important role in elasticity theory. %; and it is also called {\it strong ellipticity} in the community of elasticity theory.
In Section 4, following the ideas given in \cite{HDQ-09,QDH-09}, we further discuss positive definiteness of fourth order three dimensional (strongly) paired symmetric tensors and propose several necessary and sufficient conditions for which the concerned tensor is positive definite. Furthermore, we extend the related results to the case of sixth order three dimensional (strongly) paired symmetric tensor.

A polynomial with real coefficients is called a sum of squares (SOS) if it can be expressed as a sum of several squares of polynomials with real coefficients \cite{Laurent-09,NDS-06}. It is obvious that a polynomial is positive semidefinite if it is an SOS polynomial. In Section 5, we investigate the SOS properties of polynomials defined by a fourth order three dimensional or a sixth order three dimensional (strongly) paired symmetric tensor. We give several necessary and/or sufficient conditions of a fourth order three dimensional (strongly) paired symmetric tensor being an SOS tensor, and propose several conditions under which a fourth order three dimensional (strongly) paired symmetric tensor is positive definite. In particular, we extend the related results to the case of sixth order three dimensional (strongly) paired symmetric tensor.

In \cite{HHQ-13}, the authors introduced the tensor conic linear programming problem and proposed a sequential semidefinite programming method to solve it. As an application, they showed that their method can be applied to find the smallest $Z$-eigenvalue of a symmetric tensor. In Section 6, by using the special structure of the (strongly) paired symmetric tensor, we propose a sequential semidefinite programming method to compute the smallest $M$-eigenvalue of a fourth order three dimensional (strongly) paired symmetric tensor, which is an extension of the method proposed in \cite{HHQ-13}. By this method, we can check whether a fourth order three dimensional (strongly) paired symmetric tensor is positive definite or not.

In Section 7, we give some numerical results of our methods for judging whether a fourth order three dimensional or sixth order three dimensional (strongly) paired symmetric tensor is positive definite or not. The preliminary numerical results are consistent with our theoretical results.

Some concluding remarks are made in Section 8.

In the remaining parts of our paper, we will simply call a three dimensional strongly paired symmetric tensor an
elasticity tensor as the main motivation of our paper is the strong ellipticity of elasticity and higher order elasticity tensors.

\section{Preliminaries}

In this section, we consider fourth order three dimensional and sixth order three dimensional paired symmetric (elasticity) tensor and discuss related basic properties.

\noindent{\bf 2.1 Fourth order paired symmetric tensor}. For any $\mathscr{A}=(a_{ijkl})\in \mathbb{T}_{4,3}$, if
\begin{eqnarray*}%\label{e-4-smy-1}
a_{ijkl}=a_{jikl}=a_{ijlk}=a_{jilk},\quad \forall i,j,k,l\in \{1,2,3\},
\end{eqnarray*}
then $\mathscr{A}$ is a {\it paired symmetric tensor} $\mathscr{A}\in \mathbb{T}_{4,3}$; and if
\begin{eqnarray*}%\label{e-4-smy-2}
a_{ijkl}=a_{jikl}=a_{ijlk}=a_{jilk}\quad\mbox{\rm and}\quad a_{ijkl}=a_{klij},\quad \forall i,j,k,l\in \{1,2,3\},
\end{eqnarray*}
then $\mathscr{A}$ is a {\it fourth order three dimensional elasticity tensor} \cite{Ha-07}.

For any tensor $\mathscr{A}=(a_{ijkl})\in \mathbb{T}_{4,3}$, the corresponding biquadratic form is defined by
\begin{eqnarray}\label{e-4-poly}
\mathscr{A}x^2y^2:=\sum_{i,j,k,l=1}^3a_{ijkl}x_ix_jy_ky_l, \quad\forall x,y\in \mathbb{R}^3.
\end{eqnarray}
Define
\begin{eqnarray}\label{e-4-cone}
\mathcal{C}:=\{\mathscr{A}\in \mathbb{T}_{4,3}: \mathscr{A}x^2y^2\geq 0, \forall x,y\in \mathbb{R}^3\}.
\end{eqnarray}
Then, by a similar way as those in \cite{QY-14,QYW-10}, we can obtain the following results.
\begin{Proposition}\label{pro-4-basic}
Suppose that the polynomial $\mathscr{A}x^2y^2$ and the set $\mathcal{C}$ are defined by (\ref{e-4-poly}) and (\ref{e-4-cone}), respectively. Then, the following statements hold.
\begin{itemize}
\item[(i)] The interior of $\mathcal{C}$, denoted by int$\mathcal{C}$, is nonempty, and
\begin{eqnarray}\label{e-4-cone-1}
\mbox{\rm int}\mathcal{C}=\{\mathscr{A}\in \mathbb{T}_{4,3}: \mathscr{A}x^2y^2>0, \forall x,y\in \mathbb{R}^3\backslash\{0\}\}.
\end{eqnarray}
\item[(ii)] The set $\mathcal{C}$ is a pointed closed convex cone.
\end{itemize}
\end{Proposition}

\noindent {\bf Proof}. (i) It is obvious that $\mathcal{C}$ has nonempty interior. We show that (\ref{e-4-cone-1}) holds. On the one hand, if $\mathscr{A}\in \mathcal{C}$ is not positive definite, then there exists two nonzero vectors $x,y\in \mathbb{R}^3$ such that $\mathscr{A}x^2y^2=0$; and hence, for any $\epsilon>0$,
$$
(\mathscr{A}-\epsilon \mathscr{E})x^2y^2=-\epsilon(x^\top x)(y^\top y)<0
$$
where $\mathscr{E}=(e_{ijkl})\in \mathbb{T}_{4,3}$ is defined by
\begin{eqnarray}\label{e-partial-identity}
e_{ijkl}:=\left\{\begin{array}{ll}
1 & \mbox{\rm if}\; i=j, k=l,\\
0 & \mbox{\rm otherwise},
\end{array}\right.\quad \forall i,j,k,l\in \{1,2,3\}.
\end{eqnarray}
This implies that $\mathscr{A}\not\in \mbox{\rm int}\mathcal{C}$; and hence,  $\mbox{\rm int}\mathcal{C}\subseteq\{\mathscr{A}\in \mathbb{T}_{4,3}: \mathscr{A}x^2y^2>0, \forall x,y\in \mathbb{R}^3\backslash\{0\}\}$. On the other hand, if $\mathscr{A}\not\in \mbox{\rm int}\mathcal{C}$, then there exists a sequence $\{(\mathscr{B}^{(k)}, \epsilon_{(k)})\}$ satisfying $\|\mathscr{B}^{(k)}\|_{HS}=1$ and $\epsilon_{(k)}>0$ for all $k\in \{1,2,\ldots\}$ such that
\begin{eqnarray*}
\mathscr{A}+\epsilon_{(k)}\mathscr{B}^{(k)}\not\in \mathcal{C}\; \mbox{\rm for all}\; k\in \{1,2,\ldots\}\quad \mbox{\rm and}\quad \lim\limits_{k\rightarrow \infty}\epsilon_{(k)}=0,
\end{eqnarray*}
which leads to that there exist $\{x^{(k)}\},\{y^{(k)}\}\subseteq \mathbb{R}^3$ satisfying $\|x^{(k)}\|=1$ and $\|y^{(k)}\|=1$ for all $k\in \{1,2,\ldots\}$ such that
$$
\left(\mathscr{A}+\epsilon_{(k)}\mathscr{B}^{(k)}\right)\left(x^{(k)}\right)^2\left(y^{(k)}\right)^2\leq 0.
$$
Let $x^*,y^*$ be the limiting points of $\{x^{(k)}\}$ and $\{y^{(k)}\}$, respectively. Then, $\|x^*\|=1, \|y^*\|=1$ and $\mathscr{A}(x^*)^2(y^*)^2\leq 0$, which implies that $\mathscr{A}$ is not positive definite; and hence,  $\{\mathscr{A}\in \mathbb{T}_{4,3}: \mathscr{A}x^2y^2>0, \forall x,y\in \mathbb{R}^3\backslash\{0\}\}\subseteq\mbox{\rm int}\mathcal{C}$. So, (\ref{e-4-cone-1}) holds.

(ii) For any $\mathscr{A},\mathscr{B}\in \mathcal{C}$ and $\alpha,\beta\geq 0$, let $\mathscr{C}:=\alpha\mathscr{A}+\beta\mathscr{B}$. Then, for any $x,y\in \mathbb{R}^3$,
\begin{eqnarray*}
\mathscr{C}x^2y^2=(\alpha\mathscr{A}+\beta\mathscr{B})x^2y^2=\alpha\mathscr{A}x^2y^2+\beta\mathscr{B}x^2y^2\geq 0,
\end{eqnarray*}
which implies that $\mathcal{C}$ is a convex cone.

For any $\mathscr{A}=(a_{ijkl})\in \mathcal{C}$, if $-\mathscr{A}\in \mathcal{C}$, then it follows from (\ref{e-4-poly}) that
$$
\sum_{i,j,k,l=1}^3a_{ijkl}x_ix_jy_ky_l\equiv 0, \forall x,y\in \mathbb{R}^3,
$$
which yields that $\mathscr{A}=0$ by arbitrariness of $x$ and $y$. This implies that $\mathcal{C}$ is a pointed cone.

For any $\{\mathscr{A}^{(k)}\}\subseteq \mathcal{C}$ and $\lim_{k\rightarrow \infty}\mathscr{A}^{(k)}=\mathscr{A}$, it is easy to see that for any $x,y\in \mathbb{R}^n$,
$$
\mathscr{A}x^2y^2=\lim\limits_{k\rightarrow \infty}\mathscr{A}^{(k)}x^2y^2\geq 0,
$$
which implies that $\mathcal{C}$ is a closed cone.

Therefore, by combining (i) with (ii), we conclude that the results of proposition hold.     \ep

In the following, we define a class of matrices which is related to tensor $\mathscr{A}\in \mathbb{T}_{4,3}$.
\begin{Definition}\label{def-unfolded-m-4}
For any $\mathscr{A}=(a_{ijkl})\in \mathbb{T}_{4,3}$, we define a matrix by
\begin{eqnarray}\label{e-unfold-matrix}
M=(m_{st})\;\;\mbox{\rm with}\;\; m_{st}=a_{i_si_tj_sj_t}\; \forall s,t\in \{1,2,\ldots,9\}
\end{eqnarray}
where $i_1i_2\cdots i_9$ and $j_1j_2\cdots j_9$ are two arbitrary permutations of $123123123$. We say that the matrix $M$ defined by (\ref{e-unfold-matrix}) is an {\bf unfolded matrix} of tensor $\mathscr{A}$ with respect to indices $i_1i_2\cdots i_9$ and $j_1j_2\cdots j_9$.
\end{Definition}

Obviously, there exists a unique unfold matrix of tensor $\mathscr{A}$ for each pair of permutations of $i_1i_2\cdots i_9$ and $j_1j_2\cdots j_9$. The following result is easy to be obtained.

\begin{Proposition}\label{pro-matrix-0}
Suppose that $\mathscr{A}\in \mathbb{T}_{4,3}$ is a paired symmetric tensor, and $M$ is its unfolded matrix defined by Definition \ref{def-unfolded-m-4}. Then, the matrix $M$ is symmetric.
\end{Proposition}

For any $\mathscr{A}\in \mathbb{T}_{4,3}$, we define
\begin{itemize}
\item matrix $M^1=(m_{st})\in \mathbb{R}^{9\times 9}$ by
\begin{eqnarray}\label{e-matrix-A}
M^1:=(m^1_{st})\quad \mbox{\rm with}\; m^1_{st}=a_{3(i-1)+k,3(j-1)+l}=a_{ijkl}\; \mbox{\rm for any}\; i,j,k,l\in \{1,2,3\};
\end{eqnarray}
\item and matrix $M^2=(m_{st})\in \mathbb{R}^{9\times 9}$ by
\begin{eqnarray}\label{e-matrix-B}
M^2:=(m^2_{st})\quad \mbox{\rm with}\; m^2_{st}=a_{3(k-1)+i,3(l-1)+j}=a_{ijkl}\; \mbox{\rm for any}\; i,j,k,l\in \{1,2,3\}.
\end{eqnarray}
\end{itemize}
Then, $M^1$ and $M^2$ are two unfolded matrices of $\mathscr{A}$.

Moreover, for any $i,j,k,l\in \{1,2,3\}$, we define two block sub-matrices of tensor $\mathscr{A}$ by
\begin{eqnarray}\label{e-matrix-AB-sub}
A_{ij}:=(a_{ijkl})_{kl}\quad \mbox{\rm and}\quad B_{kl}:=(a_{ijkl})_{ij}.
\end{eqnarray}
Then, we can easily obtain the following results.
\begin{Proposition}\label{pro-matrix-1}
For any $\mathscr{A}\in \mathbb{T}_{4,3}$, suppose that matrices $M^1$, $M^2$, $A_{ij}$ and $B_{kl}$ are given by (\ref{e-matrix-A}), (\ref{e-matrix-B}) and (\ref{e-matrix-AB-sub}), respectively.
\begin{itemize}
\item[(i)] If $\mathscr{A}\in \mathbb{T}_{4,3}$ is a paired symmetric tensor, then the following statements hold.

(a) The matrices $M^1$ and $M^2$ are symmetric, and
    \begin{eqnarray}\label{e-matrix-AB-sub1}
    M^1=\left[\begin{array}{ccc} A_{11} & A_{12} & A_{13}\\ A_{21} & A_{22} & A_{23}\\ A_{31} & A_{32} & A_{33} \end{array}\right]\quad\mbox{\rm and}\quad
    M^2=\left[\begin{array}{ccc} B_{11} & B_{12} & B_{13}\\ B_{21} & B_{22} & B_{23}\\ B_{31} & B_{32} & B_{33} \end{array}\right];
    \end{eqnarray}

    (b) All sub-matrices $A_{ij}$ and $B_{kl}$ are symmetric.
\item[(ii)] If $\mathscr{A}\in \mathbb{T}_{4,3}$ is an elasticity tensor, then $M^1=M^2$.
\end{itemize}
\end{Proposition}

\noindent{\bf 2.2 Sixth order paired symmetric tensor}. $\mathscr{A}=(a_{ijklpq})\in \mathbb{T}_{6,3}$ is a {\it paired symmetric tensor} if its entries satisfy
\begin{eqnarray}\label{e-6-sym-1}
a_{ijklpq}=a_{jiklpq}=a_{ijlkpq}=a_{ijklqp},\quad \forall i,j,k,l,p,q\in \{1,2,3\}.
\end{eqnarray}
Furthermore, a paired symmetric tensor $\mathscr{A}\in \mathbb{T}_{6,3}$ is a {\it sixth order elasticity tensor} \cite{Ha-07} if
\begin{eqnarray}\label{e-6-sym-2}
a_{ijklpq}=a_{klijpq}=a_{ijpqkl}, \quad \forall i,j,k,l,p,q\in \{1,2,3\}.
\end{eqnarray}

For any tensor $\mathscr{A}=(a_{ijklpq})\in \mathbb{T}_{6,3}$, the corresponding homogeneous polynomial is defined by
\begin{eqnarray}\label{e-6-poly}
\mathscr{A}x^2y^2z^2:=\sum_{i,j,k,l,p,q=1}^3a_{ijklpq}x_ix_jy_ky_lz_pz_q,\quad \forall x,y,z\in \mathbb{R}^3.
\end{eqnarray}
Then, similar to Proposition \ref{pro-4-basic}, we have the following results.
\begin{Proposition}\label{pro-6-basic}
Suppose that $\mathscr{A}x^2y^2z^2$ is defined by (\ref{e-6-poly}). Then, the following statements hold.
\begin{itemize}
\item[(i)] The interior of set $\mathcal{D}$ defined by $\mathcal{D}:=\{\mathscr{A}\in \mathbb{T}_{6,3}: \mathscr{A}x^2y^2z^2\geq 0, \forall x,y,z\in \mathbb{R}^3\}$ is nonempty, and
$\mbox{\rm int}\mathcal{D}=\{\mathscr{A}\in \mathbb{T}_{6,3}: \mathscr{A}x^2y^2z^2>0, \forall x,y,z\in \mathbb{R}^3\backslash\{0\}\}$.
\item[(ii)] The set $\mathcal{D}$ is a pointed closed convex cone.
\end{itemize}
\end{Proposition}

Similar to Definition \ref{def-unfolded-m-4}, we define the unfolded matrix of sixth order three dimensional tensor as follows.
\begin{Definition}\label{def-unfolded-m-6}
For any $\mathscr{A}=(a_{ijkl}pq)\in \mathbb{T}_{6,3}$, we define a matrix by
\begin{eqnarray}\label{e-unfold-matrix-6}
N=(n_{st})\;\;\mbox{\rm with}\;\; n_{st}=a_{i_si_tj_sj_tk_sk_t}\; \forall s,t\in \{1,2,\ldots,27\}
\end{eqnarray}
where $i_1i_2\cdots i_{27}$, $j_1j_2\cdots j_{27}$ and  $k_1k_2\cdots k_{27}$ are three arbitrary permutations of $\underbrace{123123\cdot 123}_{27}$. We say that the matrix $N$ defined by (\ref{e-unfold-matrix-6}) is an {\bf unfolded matrix} of tensor $\mathscr{A}$ with respect to indices $i_1i_2\cdots i_{27}$, $j_1j_2\cdots j_{27}$ and  $k_1k_2\cdots k_{27}$.
\end{Definition}

Obviously, there exists a unique unfolded matrix of tensor $\mathscr{A}$ for each triple of permutations of $i_1i_2\cdots i_{27}$, $j_1j_2\cdots j_{27}$ and  $k_1k_2\cdots k_{27}$; and the following result holds.

\begin{Proposition}\label{pro-matrix-add0}
Suppose that $\mathscr{A}\in \mathbb{T}_{6,3}$ is a paired symmetric tensor, and $N$ is its unfolded matrix defined by Definition \ref{def-unfolded-m-6}. Then, the matrix $N$ is symmetric.
\end{Proposition}

In the following, we give several specific examples of the unfolded matrix. That is, for any $\mathscr{A}\in \mathbb{T}_{6,3}$, we define
\begin{eqnarray}\label{e-6-matrix-1}
\begin{array}{rcl}
N^1&:=&(n_{st})\in \mathbb{R}^{27\times 27}\; \mbox{\rm with}\\
 &  & n_{st}=n_{3[3(i-1)+(k-1)]+p,3[3(j-1)+(l-1)]+q}=a_{ijklpq},\; \forall i,j,k,l,p,q\in \{1,2,3\};\\
N^2&:=&(n_{st})\in \mathbb{R}^{27\times 27}\; \mbox{\rm with}\\
 &  & n_{st}=n_{3[3(k-1)+(i-1)]+p,3[3(l-1)+(j-1)]+q}=a_{ijklpq},\; \forall i,j,k,l,p,q\in \{1,2,3\};\\
N^3&:=&(n_{st})\in \mathbb{R}^{27\times 27}\; \mbox{\rm with}\\
 &  & n_{st}=n_{3[3(p-1)+(k-1)]+i,3[3(q-1)+(l-1)]+j}=a_{ijklpq},\; \forall i,j,k,l,p,q\in \{1,2,3\};\\
N^4&:=&(n_{st})\in \mathbb{R}^{27\times 27}\; \mbox{\rm with}\\
 &  & n_{st}=n_{3[3(i-1)+(p-1)]+k,3[3(j-1)+(q-1)]+l}=a_{ijklpq},\; \forall i,j,k,l,p,q\in \{1,2,3\};\\
N^5&:=&(n_{st})\in \mathbb{R}^{27\times 27}\; \mbox{\rm with}\\
 &  & n_{st}=n_{3[3(k-1)+(p-1)]+i,3[3(l-1)+(q-1)+j}=a_{ijklpq},\; \forall i,j,k,l,p,q\in \{1,2,3\};\\
N^6&:=&(n_{st})\in \mathbb{R}^{27\times 27}\; \mbox{\rm with}\\
 &  & n_{st}=f_{3[3(p-1)+(i-1)]+k,3[3(q-1)+(j-1)]+l}=a_{ijklpq},\; \forall i,j,k,l,p,q\in \{1,2,3\}.
\end{array}
\end{eqnarray}
Then, we have the following results.
\begin{Proposition}\label{pro-6-matrix-1}
For any paired symmetric tensor $\mathscr{A}\in \mathbb{T}_{6,3}$, we have the following results.
\begin{itemize}
\item All matrices $N^1, N^2, \ldots, N^6$ defined by (\ref{e-6-matrix-1}) are symmetric.
\item If $\mathscr{A}$ is an elasticity tensor, then $N^1=N^2=\cdots=N^6$.
\end{itemize}
\end{Proposition}

For any tensor $\mathscr{A}=(a_{ijklpq})\in \mathbb{T}_{6,3}$ and any $i,j,k,l,p,q\in \{1,2,3\}$, we define three block sub-tensors of tensor $\mathscr{A}$ by
\begin{eqnarray}\label{e-six-subtensor-d}
\mathscr{A}_{ij}:=(a_{ijklpq})_{klpq},\quad \mathscr{B}_{kl}:=(a_{ijklpq})_{ijpq}\quad \mbox{\rm and}\quad \mathscr{C}_{pq}:=(a_{ijklpq})_{ijkl}.
\end{eqnarray}
Then, we have the following results.

\begin{Proposition}\label{pro-six-tensor-1}
For any $\mathscr{A}\in \mathbb{T}_{6,3}$, let sub-tensors $\mathscr{A}_{ij}$, $\mathscr{B}_{kl}$ and $\mathscr{C}_{pq}$ be defined by (\ref{e-six-subtensor-d}).
\begin{itemize}
\item If $\mathscr{A}$ is a paired symmetric tensor, then all sub-tensors $\mathscr{A}_{ij}$, $\mathscr{B}_{kl}$ and $\mathscr{C}_{pq}$ are paired symmetric tensors.
\item If $\mathscr{A}$ is an elasticity tensor (i.e., its entries satisfy (\ref{e-6-sym-1}) and (\ref{e-6-sym-2})), then $\mathscr{A}_{st}=\mathscr{B}_{st}=\mathscr{C}_{st}$ for all $s,t\in \{1,2,3\}$.
\end{itemize}
\end{Proposition}

\section{$M$-Eigenvalue and Properties}

In this section, we extend the concept of $M$-eigenvalues for fourth order paired symmetric tensor introduced in \cite{HDQ-09,QDH-09} to sixth order three dimensional paired symmetric (elasticity) tensor and bi-block symmetric tensor, and discuss some related properties.

For any paired symmetric tensor $\mathscr{A}=(a_{ijklpq})\in \mathbb{T}_{6,3}$, the corresponding homogeneous polynomial is given in
(\ref{e-6-poly}). For any $x,y,z\in \mathbb{R}^3$, we let $\mathscr{A}xy^2z^2,\mathscr{A}x^2yz^2,\mathscr{A}x^2y^2z\in \mathbb{R}^3$ be defined by
\begin{eqnarray*}
%\begin{array}{rcl}
(\mathscr{A}xy^2z^2)_i&:=&\sum_{j,k,l,p,q=1}^3a_{ijklpq}x_jy_ky_lz_pz_q,\quad \forall i\in \{1,2,3\},\\
(\mathscr{A}x^2yz^2)_k&:=&\sum_{i,j,l,p,q=1}^3a_{ijklpq}x_ix_jy_lz_pz_q,\quad \forall k\in \{1,2,3\},\\
(\mathscr{A}x^2y^2z)_p&:=& \sum_{i,j,k,l,q=1}^3a_{ijklpq}x_ix_jy_ky_lz_q,\quad \forall p\in \{1,2,3\}.
%\end{array}
\end{eqnarray*}
Then, it is easy to see that
\begin{eqnarray}\label{e-property-1}
\langle x,\mathscr{A}xy^2z^2\rangle=\mathscr{A}x^2y^2z^2,\quad \langle y,\mathscr{A}x^2yz^2\rangle=\mathscr{A}x^2y^2z^2,\quad \langle z,\mathscr{A}x^2y^2z\rangle=\mathscr{A}x^2y^2z^2
\end{eqnarray}
hold for any $x,y,z\in \mathbb{R}^3$.

\begin{Definition}\label{def-6-M-eigen}
For any paired symmetric tensor $\mathscr{A}=(a_{ijklpq})\in \mathbb{T}_{6,3}$, if there exist $\lambda\in \mathbb{R}$ and $x,y,z\in \mathbb{R}^3$ such that
$$
\left\{\begin{array}{l}
\mathscr{A}xy^2z^2=\lambda x,\quad \mathscr{A}x^2yz^2=\lambda y,\quad \mathscr{A}x^2y^2z=\lambda z,\\
x^\top x=1,\quad y^\top y=1,\quad z^\top z=1,
\end{array}\right.
$$
then $\lambda$ is called an $M$-eigenvalue of $\mathscr{A}$ and $x,y,z$ are the eigenvectors of $\mathscr{A}$ associated with the $M$-eigenvalue $\lambda$.
\end{Definition}

\begin{Theorem}\label{thm-Meig-exist}
For any paired symmetric (elasticity) tensor $\mathscr{A}=(a_{ijklpq})\in \mathbb{T}_{6,3}$, its $M$-eigenvalues always exist. Moreover, if $x,y,z$ are the eigenvectors of $\mathscr{A}$ associated with the $M$-eigenvalue $\lambda$, then $\lambda=\mathscr{A}x^2y^2z^2$.
\end{Theorem}

\noindent {\bf Proof}. Let $\mathscr{A}x^2y^2z^2$ be defined by (\ref{e-6-poly}). We consider the following optimization problem:
\begin{eqnarray}\label{e-6-opt-p}
\left\{\begin{array}{cl}
\min & \mathscr{A}x^2y^2z^2\\
\mbox{\rm s.t.} & x^\top x=1,\quad y^\top y=1,\quad z^\top z=1.
\end{array}\right.
\end{eqnarray}
It is easy to see that the feasible set of (\ref{e-6-opt-p}) is compact and the objective function of (\ref{e-6-opt-p}) is continuous. Thus, the optimization problem (\ref{e-6-opt-p}) has at least a minimizer, say $(x^*,y^*,z^*)$, which satisfies the first order optimality condition of (\ref{e-6-opt-p}), i.e., there exist $\alpha, \beta,\gamma\in \mathbb{R}$ such that
$$
\left\{\begin{array}{l}
\mathscr{A}x^*(y^*)^2(z^*)^2=\alpha x^*,\quad \mathscr{A}(x^*)^2y^*(z^*)^2=\beta y^*,\quad \mathscr{A}(x^*)^2(y^*)^2z^*=\gamma z^*,\\
(x^*)^\top x^*=1,\quad (y^*)^\top y^*=1,\quad (z^*)^\top z^*=1.
\end{array}\right.
$$
This, together with (\ref{e-property-1}), implies that
$$
\alpha=\beta=\gamma=\mathscr{A}(x^*)^2(y^*)^2(z^*)^2.
$$
Thus, $\alpha$ is an $M$-eigenvalue of $\mathscr{A}$ and $x^*,y^*,z^*$ are the eigenvectors of $\mathscr{A}$ associated with the $M$-eigenvalue $\alpha$. We complete the proof. \ep

By Theorem \ref{thm-Meig-exist}, we have the following result.

\begin{Theorem}\label{thm-6-positive-basic}
A paired symmetric (elasticity) tensor $\mathscr{A}=(a_{ijklpq})\in \mathbb{T}_{6,3}$ is positive definite if and only if the smallest $M$-eigenvalue of $\mathscr{A}$ is positive.
\end{Theorem}

This theorem demonstrates that positive definiteness detection of a paired symmetric (elasticity) tensor $\mathscr{A}\in \mathbb{T}_{6,3}$ can be done by computing the smallest $M$-eigenvalue of $\mathscr{A}$. Moreover, from Definition \ref{def-6-M-eigen} and Theorem \ref{thm-6-positive-basic}, it is easy to obtain the following results.

\begin{Theorem}\label{thm-6-positive-basic1}
For any paired symmetric (elasticity) tensor $\mathscr{A}=(a_{ijklpq})\in \mathbb{T}_{6,3}$, it follows that $\lambda$ is an $M$-eigenvalue of $\mathscr{A}$ if and only if $-\lambda$ is an $M$-eigenvalue of $-\mathscr{A}$; and furthermore, $\mathscr{A}$ is positive definite if and only if the largest $M$-eigenvalue of $-\mathscr{A}$ is negative.
\end{Theorem}

From the point of view of numerical calculation, Theorem \ref{thm-6-positive-basic1} is useful since positive definiteness detection of a paired symmetric (elasticity) tensor $\mathscr{A}\in \mathbb{T}_{6,3}$ can be done by computing the largest $M$-eigenvalue of $-\mathscr{A}$.

For any $(\lambda,x,y,z)\in \mathbb{R}\times \mathbb{R}^3\times \mathbb{R}^3\times \mathbb{R}^3$ and tensor $\mathscr{A}\in \mathbb{T}_{6,3}$, $\lambda x^2y^2z^2$ is a rank-one sixth order paired symmetric tensor with its entries being $\lambda x_ix_jy_ky_lz_pz_q$ for all $i,j,k,l,p,q\in \{1,2,3\}$. We say that $\lambda_* (x^*)^2(y^*)^2(z^*)^2$ is the best rank-one approximation of $\mathscr{A}$ if
$(\lambda_*,x^*,y^*,z^*)\in \mathbb{R}\times \mathbb{R}^3\times \mathbb{R}^3\times \mathbb{R}^3$ solves the optimization problem:
\begin{eqnarray}\label{e-best-1-appr}
\begin{array}{cl}
\min & \|\mathscr{A}-\lambda x^2y^2z^2\|_{HS}^2\\
\mbox{\rm s.t.} & \lambda\in \mathbb{R} \quad\mbox{\rm and}\quad x^\top x=1,\; y^\top y=1,\; z^\top z=1,\; \forall x,y,z\in \mathbb{R}^3.
\end{array}
\end{eqnarray}
The best rank-one approximation has wide applications in signal and image processing, wireless communication systems, independent component analysis, and so on.

\begin{Theorem}\label{thm-best-1-appr}
For any paired symmetric (elasticity) tensor $\mathscr{A}=(a_{ijklpq})\in \mathbb{T}_{6,3}$, if $\lambda_*$ is an $M$-eigenvalue of $\mathscr{A}$ with the largest absolute value among all $M$-eigenvalues of $\mathscr{A}$, and $x^*,y^*,z^*\in \mathbb{R}^3$ are the eigenvectors of $\mathscr{A}$ associated with the $M$-eigenvalue $\lambda_*$, then $\lambda_*(x^*)^2(y^*)^2(z^*)^2$ is the best rank-one approximation of $\mathscr{A}$.
\end{Theorem}

\noindent {\bf Proof}. We denote the feasible set of (\ref{e-best-1-appr}) by $\Omega$, i.e.,
$$
\Omega:=\{(\lambda,x,y,z)\in \mathbb{R}\times \mathbb{R}^3\times \mathbb{R}^3\times \mathbb{R}^3:\; x^\top x=1,\; y^\top y=1,\; z^\top z=1\}.
$$
On the one hand, since
\begin{eqnarray*}
\begin{array}{l}
\min \{\|\mathscr{A}-\lambda x^2y^2z^2\|_{HS}^2:\; (\lambda,x,y,z)\in \Omega\}\\
\quad =\min\{\|\mathscr{A}\|_{HS}^2-2\lambda\mathscr{A}x^2y^2z^2+\lambda^2(x^\top x)(y^\top y)(z^\top z):\; (\lambda,x,y,z)\in \Omega\}\\
\quad =\min\{\|\mathscr{A}\|_{HS}^2-2\lambda\mathscr{A}x^2y^2z^2+\lambda^2:\; (\lambda,x,y,z)\in \Omega\}
\end{array}
\end{eqnarray*}
and when $\lambda=\mathscr{A}xxyyzz$,
\begin{eqnarray*}
\begin{array}{rcl}
\min \{\|\mathscr{A}-\lambda x^2y^2z^2\|_{HS}^2:\; (\lambda,x,y,z)\in \Omega\}
&=&\min\{\|\mathscr{A}\|_{HS}^2-(\lambda\mathscr{A}x^2y^2z^2)^2:\; (\lambda,x,y,z)\in \Omega\}\\
&=&\|\mathscr{A}\|_{HS}^2-\max\{(\mathscr{A}x^2y^2z^2)^2:\; (\lambda,x,y,z)\in \Omega\},
\end{array}
\end{eqnarray*}
it is easy to see that $\|\mathscr{A}-\lambda x^2y^2z^2\|_{HS}^2$ is the smallest if and only if the absolute value of $\lambda$ is the largest among all $\lambda\in \mathbb{R}$ satisfying $\lambda=\mathscr{A}x^2y^2z^2$ and $(\lambda,x,y,z)\in \Omega$. On the other hand, since $\lambda_*$ is an $M$-eigenvalues of $\mathscr{A}$ and $x^*,y^*,z^*\in \mathbb{R}^3$ are the eigenvectors of $\mathscr{A}$ associated with the $M$-eigenvalue $\lambda_*$, from Definition \ref{def-6-M-eigen} it follows that $(\lambda_*,x^*,y^*,z^*)\in \Omega$; and from Theorem \ref{thm-Meig-exist} it follows that $\lambda_*=\mathscr{A}(x^*)^2(y^*)^2(z^*)^2$.

Therefore, by the assumption that $\lambda_*$ is an $M$-eigenvalues of $\mathscr{A}$ with the largest absolute value among all $M$-eigenvalues of $\mathscr{A}$, we obtain that $\lambda_*(x^*)^2(y^*)^2(z^*)^2$ is the best rank-one approximation of $\mathscr{A}$, which completes the proof.  \ep

In the following, we introduce the {\it bi-block symmetric tensor} and study some related properties.
\begin{Definition}
$\mathscr{A}=(a_{i_1i_2\cdots i_{2m}})\in \mathbb{T}_{2m,3}$ is called a bi-block symmetric tensor if its indices $\{i_1,i_2,\cdots, i_{2m}\}$ are divided into two adjacent blocks $\{i_1,i_2,\cdots, i_t\}$ and $\{i_{t+1},i_{t+2},\cdots, i_{2m}\}$ with $t\in [1, 2m]$ being an even number and entries of $\mathscr{A}$ being invariant under any permutation of indices in every block of $\{i_1,i_2,\cdots, i_{t}\}$ and $\{i_{t+1},i_{t+2},\cdots, i_{2m}\}$, i.e.,
\begin{eqnarray}\label{e-t-sym-1}
a_{i_1i_2\cdots i_{t}i_{t+1}i_{t+2}\cdots i_{2m}}=a_{\sigma(i_1i_2\cdots i_{t})\sigma(i_{t+1}i_{t+2}\cdots i_{2m})}
\end{eqnarray}
for all $i_1,i_2,\ldots,i_{2m}\in \{1,2,3\}$, where $\sigma(i_1i_2\cdots i_{t})$ denotes an arbitrary permutation of $i_1i_2\cdots i_{t}$.
\end{Definition}

For any even number $t\in [1,2m]$ and $x,y\in \mathbb{R}^3$, we use the following notation:
\begin{eqnarray*}
\begin{array}{l}
\qquad x^{[t]}:=(x_1^t,x_2^t,x_3^t)^\top;\\
\qquad \mathscr{A}x^ty^{2m-t}:=\sum\limits_{i_1,\cdots, i_t, i_{t+1},\cdots,i_{2m}=1}^3a_{i_1\cdots i_ti_{t+1}\cdots i_{2m}}x_{i_1}\cdots x_{t}y_{t+1}\cdots y_{2m}\\
\qquad \mathscr{A}x^{t-1}y^{2m-t}\in \mathbb{R}^3\; \mbox{\rm with}\\
\qquad\quad (\mathscr{A}x^{t-1}y^{2m-t})_i:=\sum\limits_{i_2,\cdots, i_t, i_{t+1},\cdots,i_{2m}=1}^3a_{ii_2\cdots i_ti_{t+1}\cdots i_{2m}}x_{i_2}\cdots x_{t}y_{t+1}\cdots y_{2m},\; \forall i\in \{1,2,3\},\\
\qquad \mathscr{A}x^ty^{2m-t-1}\in \mathbb{R}^3\; \mbox{\rm with}\\
\qquad\quad  (\mathscr{A}x^ty^{2m-t-1})_i:=\sum\limits_{i_1,\cdots, i_t, ii_{t+2},\cdots,i_{2m}=1}^3a_{i_1\cdots i_tii_{t+2}\cdots i_{2m}}x_{i_1}\cdots x_{t}y_{t+2}\cdots y_{2m},\; \forall i\in \{1,2,3\}.
\end{array}
\end{eqnarray*}

\begin{Definition}\label{def-t-M-eigen}
For any bi-block symmetric tensor $\mathscr{A}=(a_{i_1\cdots i_ti_{t+1}\cdots i_{2m}})\in \mathbb{T}_{2m,3}$ whose entries satisfy bi-block symmetry given by (\ref{e-t-sym-1}), if there exist $\lambda\in \mathbb{R}$ and $x,y\in \mathbb{R}^3$ such that
$$
\left\{\begin{array}{l}
\mathscr{A}x^{t-1}y^{2m-t}=\lambda x,\quad \mathscr{A}x^ty^{2m-t-1}=\lambda y,\\
\left(x^{\left[\frac{t}{2}\right]}\right)^\top x^{\left[\frac{t}{2}\right]}=1,\quad \left(y^{\left[\frac{2m-t}{2}\right]}\right)^\top y^{\left[\frac{2m-t}{2}\right]}=1,
\end{array}\right.
$$
then $\lambda$ is called an $M$-eigenvalue of $\mathscr{A}$ and $x,y$ are the eigenvectors of $\mathscr{A}$ associated with the $M$-eigenvalue $\lambda$.
\end{Definition}
Then, similar to Theorem \ref{thm-Meig-exist}, we can obtain the following results.

\begin{Theorem}\label{thm-t-Meig-exist}
For any bi-block symmetric tensor $\mathscr{A}=(a_{i_1\cdots i_ti_{t+1}\cdots i_{2m}})\in \mathbb{T}_{2m,3}$ whose entries satisfy bi-block symmetry given by (\ref{e-t-sym-1}), $M$-eigenvalues always exist. Moreover, if $x,y$ are the eigenvectors of $\mathscr{A}$ associated with the $M$-eigenvalue $\lambda$, then $\lambda=\mathscr{A}x^ty^{2m-t}$.
\end{Theorem}

By using Theorem \ref{thm-t-Meig-exist}, the following result holds.

\begin{Theorem}\label{thm-t-Meig-add}
A bi-block symmetric tensor $\mathscr{A}=(a_{i_1\cdots i_ti_{t+1}\cdots i_{2m}})\in \mathbb{T}_{2m,3}$ whose entries satisfy bi-block symmetry given by (\ref{e-t-sym-1}) is positive definite if and only if the smallest $M$-eigenvalue of $\mathscr{A}$ is positive.
\end{Theorem}

From Definition \ref{def-t-M-eigen} and Theorem \ref{thm-t-Meig-add}, we have the following results.

\begin{Theorem}\label{thm-t-Meig-add1}
For any bi-block symmetric tensor $\mathscr{A}\in \mathbb{T}_{2m,3}$  whose entries satisfy bi-block symmetry given by (\ref{e-t-sym-1}), it follows that $\lambda$ is an $M$-eigenvalue of $\mathscr{A}$ if and only if $-\lambda$ is an $M$-eigenvalue of $-\mathscr{A}$; and furthermore, $\mathscr{A}$ is positive definite if and only if the largest $M$-eigenvalue of $-\mathscr{A}$ is negative.
\end{Theorem}

Moreover, similar to Theorem \ref{thm-best-1-appr}, we can obtain the following result.

\begin{Theorem}\label{thm-t-best-1-appr}
For any bi-block symmetric tensor $\mathscr{A}=(a_{i_1\cdots i_ti_{t+1}\cdots i_{2m}})\in \mathbb{T}_{2m,3}$ whose entries satisfy bi-block symmetry given by (\ref{e-t-sym-1}), if $\lambda_*$ is an $M$-eigenvalues of $\mathscr{A}$ with the largest absolute value among all $M$-eigenvalues of $\mathscr{A}$, and $x^*,y^*\in \mathbb{R}^3$ are the eigenvectors of $\mathscr{A}$ associated with the $M$-eigenvalue $\lambda_*$, then $\lambda_*(x^*)^t(y^*)^{2m-t}$ is the best rank-one approximation of $\mathscr{A}$.
\end{Theorem}

For any bi-block symmetric tensor $\mathscr{A}=(a_{i_1\cdots i_ti_{t+1}\cdots i_{2m}})\in \mathbb{T}_{2m,3}$ whose entries satisfy bi-block symmetry given by (\ref{e-t-sym-1}), when $m=2$ and $t=2$, the tensor $\mathscr{A}$ reduces to a fourth order three dimensional paired symmetric tensor; and hence, the definition of $M$-eigenvalue and results of Theorems \ref{thm-t-Meig-exist}, \ref{thm-t-Meig-add} and \ref{thm-t-best-1-appr} reduce to those given in \cite{QDH-09}. For example, by Theorem \ref{thm-t-Meig-add} we have

\begin{Corollary}\label{thm-4-positive-basic}
A paired symmetric (elasticity) tensor $\mathscr{A}=(a_{ijkl})\in \mathbb{T}_{4,3}$ is positive definite if and only if the smallest $M$-eigenvalue of $\mathscr{A}$ is positive.
\end{Corollary}

\section{Eigenvalue and Positive Definiteness}

In this section, following the ideas given in \cite{HDQ-09}, we consider positive definiteness of fourth order three dimensional and sixth order three dimensional paired symmetric (elasticity) tensors. Specially, we propose several necessary and sufficient conditions under which the concerned tensors are positive definite.

\noindent{\bf 4.1 Fourth order paired symmetric tensors}. For any $x, y\in \mathbb{R}^3$, we define two matrices by
\begin{eqnarray}\label{e-biquad-twoM}
A(y):=\left[\begin{array}{ccc} y^\top A_{11}y & y^\top A_{12}y & y^\top A_{13}y \\ y^\top A_{21}y & y^\top A_{22}y & y^\top A_{23}y\\ y^\top A_{31}y & y^\top A_{32}y & y^\top A_{33}y \end{array}\right]
\quad\mbox{\rm and}\quad
B(x):=\left[\begin{array}{ccc} x^\top B_{11}x & x^\top B_{12}x & x^\top B_{13}x \\ x^\top B_{21}x & x^\top B_{22}x & x^\top B_{23}x\\ x^\top B_{31}x & x^\top B_{32}x & x^\top B_{33}x \end{array}\right],
\end{eqnarray}
where $A_{ij}$ and $B_{kl}$ for $i,j,k,l\in \{1,2,3\}$ are defined by (\ref{e-matrix-AB-sub}). Recall that for any tensor $\mathscr{A}=(a_{ijkl})\in \mathbb{T}_{4,3}$, the corresponding biquadratic form is given by (\ref{e-4-poly}), i.e.,
\begin{eqnarray*}%\label{e-biquad-form}
\mathscr{A}x^2y^2=\sum_{i,j,k,l=1}^3a_{ijkl}x_ix_jy_ky_l,\;\; \forall x,y\in \mathbb{R}^3.
\end{eqnarray*}
Thus, it follows that
\begin{eqnarray*}%\label{e-biquad-twoM-1}
\mathscr{A}x^2y^2=x^\top A(y)x=y^\top B(x)y,\;\; \forall x, y\in \mathbb{R}^3.
\end{eqnarray*}

\begin{Proposition}\label{pro-matrix-2}
For any paired symmetric (elasticity) tensor $\mathscr{A}\in \mathbb{T}_{4,3}$, suppose that the matrices $A(\cdot)$ and $B(\cdot)$ are defined by (\ref{e-biquad-twoM}). Then, the following results are equivalent.
\begin{itemize}
\item[(i)] The polynomial $\mathscr{A}x^2y^2$ defined by (\ref{e-4-poly}) is positive definite.
\item[(ii)] The matrix $A(y)$ is positive definite for all $y\in \mathbb{R}^3\backslash\{0\}$.
\item[(iii)] The matrix $B(x)$ is positive definite for all $x\in \mathbb{R}^3\backslash\{0\}$.
\end{itemize}
\end{Proposition}

It is well known that a symmetric matrix $M$ is positive definite if and only if all leading principal minors of $M$ are positive, which is the Sylvester's criterion. Thus, we have the following results, where det$(M)$ denotes the determinant of the matrix $M$.
\begin{Theorem}\label{thm-4order-1}
Suppose that $\mathscr{A}\in \mathbb{T}_{4,3}$ is a paired symmetric tensor, the polynomial $\mathscr{A}x^2y^2$ is defined by (\ref{e-4-poly}), and matrices $A(\cdot)$ and $B(\cdot)$ are defined by (\ref{e-biquad-twoM}). Then, the polynomial $\mathscr{A}x^2y^2$ is positive definite if and only if one of the following results holds.
\begin{itemize}
\item[(i)] The matrix $A_{11}$ is positive definite, and
$$
y^\top A_{11}y y^\top A_{22}y - y^\top A_{12}y y^\top A_{21}y>0\quad\mbox{\rm and}\quad \mbox{\rm det}(A(y))>0,\;\; \forall y\in \mathbb{R}^3\backslash\{0\}.
$$
\item[(ii)] The matrix $B_{11}$ is positive definite, and
$$
x^\top B_{11}x x^\top B_{22}x - x^\top B_{12}x x^\top B_{21}x>0\quad\mbox{\rm and}\quad \mbox{\rm det}(B(x))>0,\;\; \forall x\in \mathbb{R}^3\backslash\{0\}.
$$
\end{itemize}
Furthermore, if $\mathscr{A}$ is an elasticity tensor, then the above (i) and (ii) are the same.
\end{Theorem}
It is easy to see that
\begin{itemize}
\item[(a)] $y^\top A_{11}y y^\top A_{22}y - y^\top A_{12}y y^\top A_{21}y$ is a special homogeneous polynomial of degree 4, and there exists a unique symmetric tensor $\mathscr{T}_A^1\in \mathbb{T}_{4,3}$ such that $\mathscr{T}_A^1y^4=y^\top A_{11}y y^\top A_{22}y - y^\top A_{12}y y^\top A_{21}y$; and
\item[(b)] $\mbox{\rm det}(A(y))$ is a special homogeneous polynomial of degree 6, and there exists a unique symmetric tensor $\mathscr{T}_A^2\in \mathbb{T}_{6,3}$ such that $\mathscr{T}_A^2y^6=\mbox{\rm det}(A(y))$.
\end{itemize}
Similarly, there exist symmetric tensors $\mathscr{T}_B^1\in \mathbb{T}_{4,3}$ and $\mathscr{T}_B^2\in \mathbb{T}_{6,3}$ such that
\begin{itemize}
\item[(c)] $\mathscr{T}_B^1x^4=x^\top B_{11}x x^\top B_{22}x - x^\top B_{12}x x^\top B_{21}x$ and $\mathscr{T}_B^2y^6=\mbox{\rm det}(B(x))$.
\end{itemize}
By combining the theory of $Z$-eigenvalues of symmetric tensors with Theorem \ref{thm-4order-1}, we have the following results.
\begin{Theorem}\label{thm-4order-2}
Suppose that $\mathscr{A}\in \mathbb{T}_{4,3}$ is a paired symmetric tensor, the polynomial $\mathscr{A}x^2y^2$ is defined by (\ref{e-4-poly}), and symmetric tensors $\mathscr{T}_A^1,\mathscr{T}_A^2,\mathscr{T}_B^1,\mathscr{T}_B^2$ are defined by the above (a)-(c). Then, the polynomial $\mathscr{A}x^2y^2$ is positive definite if and only if one of the following results holds.
\begin{itemize}
\item[(i)] The symmetric matrix $A_{11}$ is positive definite, and the smallest $Z$-eigenvalues of the symmetric tensors $\mathscr{T}_A^1$ and $\mathscr{T}_A^2$ are positive.
\item[(ii)] The symmetric matrix $B_{11}$ is positive definite, and the smallest $Z$-eigenvalues of the symmetric tensors $\mathscr{T}_B^1$ and $\mathscr{T}_B^2$ are positive.
\end{itemize}
Furthermore, if $\mathscr{A}$ is an elasticity tensor, then the above (i) and (ii) are the same.
\end{Theorem}

\noindent{\bf 4.2 Sixth order paired symmetric tensors}. In this part, we consider positive definiteness of sixth order three dimensional paired symmetric tensors and elasticity tensors.

For any $x,y,z\in \mathbb{R}^3$, we define three matrices by
\begin{eqnarray}\label{e-biquad-threeM}
\begin{array}{rcl}
A(y,z)&:=&\left(\begin{array}{ccc} \mathscr{A}_{11}y^2z^2 & \mathscr{A}_{12}y^2z^2 & \mathscr{A}_{13}y^2z^2 \\ \mathscr{A}_{21}y^2z^2 & \mathscr{A}_{22}y^2z^2 & \mathscr{A}_{23}y^2z^2\\ \mathscr{A}_{31}y^2z^2 & \mathscr{A}_{32}y^2z^2 & \mathscr{A}_{33}y^2z^2 \end{array}\right),\\
B(x,z)&:=&\left(\begin{array}{ccc} \mathscr{B}_{11}x^2z^2 & \mathscr{B}_{12}x^2z^2 & \mathscr{B}_{13}x^2z^2 \\ \mathscr{B}_{21}x^2z^2 & \mathscr{B}_{22}x^2z^2 & \mathscr{B}_{23}x^2z^2\\ \mathscr{B}_{31}x^2z^2 & \mathscr{B}_{32}x^2z^2 & \mathscr{B}_{33}x^2z^2 \end{array}\right),\\
C(x,y)&:=&\left(\begin{array}{ccc} \mathscr{C}_{11}x^2y^2 & \mathscr{C}_{12}x^2y^2 & \mathscr{C}_{13}x^2y^2 \\ \mathscr{C}_{21}x^2y^2 & \mathscr{C}_{22}x^2y^2 & \mathscr{C}_{23}x^2y^2\\ \mathscr{C}_{31}x^2y^2 & \mathscr{C}_{32}x^2y^2 & \mathscr{C}_{33}x^2y^2 \end{array}\right).
\end{array}
\end{eqnarray}
Then, it is easy to see that
\begin{eqnarray}\label{e-six-poly-1}
\mathscr{A}x^2y^2z^2=x^\top A(y,z)x=y^\top B(x,z)y=z^\top C(x,y)z,\;\; \forall x, y,z\in \mathbb{R}^3.
\end{eqnarray}

\begin{Proposition}\label{pro-six-poly-2}
For any paired symmetric (elasticity) tensor $\mathscr{A}\in \mathbb{T}_{6,3}$, suppose that matrices $A(\cdot,\cdot)$, $B(\cdot,\cdot)$ and $C(\cdot,\cdot)$ are defined by (\ref{e-biquad-threeM}). Then, the following results are equivalent.
\begin{itemize}
\item[(i)] The polynomial $\mathscr{A}x^2y^2z^2$ defined by (\ref{e-6-poly}) is positive definite.
\item[(ii)] The matrix $A(y,z)$ is positive definite for all $y,z\in \mathbb{R}^3\backslash\{0\}$.
\item[(iii)] The matrix $B(x,z)$ is positive definite for all $x,z\in \mathbb{R}^3\backslash\{0\}$.
\item[(iv)] The matrix $C(x,y)$ is positive definite for all $x,y\in \mathbb{R}^3\backslash\{0\}$.
\end{itemize}
Furthermore, if $\mathscr{A}$ is an elasticity tensor, then the above (ii), (iii) and (iv) are the same.
\end{Proposition}

Furthermore, by Sylvester's criterion we have the following results.
\begin{Theorem}\label{thm-six-poly-1}
For any paired symmetric tensor $\mathscr{A}\in \mathbb{T}_{6,3}$, we assume that matrices $A(\cdot,\cdot)$, $B(\cdot,\cdot)$ and $C(\cdot,\cdot)$ are defined by (\ref{e-biquad-threeM}). Then, the polynomial $\mathscr{A}x^2y^2z^2$ is positive definite if and only if one of the following results holds.
\begin{itemize}
\item[(i)] $\mathscr{A}_{11}y^2z^2>0$, $(\mathscr{A}_{11}y^2z^2)(\mathscr{A}_{22}y^2z^2)-(\mathscr{A}_{12}y^2z^2)(\mathscr{A}_{21}y^2z^2)>0$ and $\mbox{\rm det}(A(y,z))>0$ for all $y,z\in \mathbb{R}^3\backslash\{0\}$.
\item[(ii)] $\mathscr{B}_{11}x^2z^2>0$, $(\mathscr{B}_{11}x^2z^2)(\mathscr{B}_{22}x^2z^2)-(\mathscr{B}_{12}x^2z^2)(\mathscr{B}_{21}x^2z^2)>0$ and $\mbox{\rm det}(B(x,z))>0$ for all $x,z\in \mathbb{R}^3\backslash\{0\}$.
\item[(iii)] $\mathscr{C}_{11}x^2y^2>0$, $(\mathscr{C}_{11}x^2y^2)(\mathscr{C}_{22}x^2y^2)-(\mathscr{C}_{12}x^2y^2)(\mathscr{C}_{21}x^2y^2)>0$ and $\mbox{\rm det}(C(x,y))>0$ for all $x,y\in \mathbb{R}^3\backslash\{0\}$.
\end{itemize}
Furthermore, if $\mathscr{A}$ is an elasticity tensor, then the above (i), (ii) and (iii) are the same.
\end{Theorem}
We consider the above result (i).
\begin{itemize}
\item It is easy to see that $\mathscr{A}_{11}$ is a paired symmetric tensor; and hence, the positive definiteness of $\mathscr{A}_{11}y^2z^2$ can be checked by the minimum $M$-eigenvalue of the tensor $\mathscr{A}_{11}$ given in the above subsection.
\item It is easy to see that $(\mathscr{A}_{11}y^2z^2)(\mathscr{A}_{22}y^2z^2)-(\mathscr{A}_{12}y^2z^2)(\mathscr{A}_{21}y^2z^2)$ is a homogeneous polynomial of degree 8 with special structure. We can define the unique paired symmetric tensor $\mathscr{T}_A^1\in \mathbb{T}_{8,3}$ such that
    $$\mathscr{T}_A^1y^4z^4=(\mathscr{A}_{11}y^2z^2)(\mathscr{A}_{22}y^2z^2)-(\mathscr{A}_{12}y^2z^2)(\mathscr{A}_{21}y^2z^2),$$
    where $(\mathscr{T}_A^1)_{i_1i_2\cdots i_8}=(\mathscr{T}_A^1)_{\sigma(i_1i_2i_3i_4)\sigma(i_5i_6i_7i_8)}$ with $\sigma(i_1i_2i_3i_4)$ being an arbitrary permutation of $i_1i_2i_3i_4$. Thus, $(\mathscr{A}_{11}y^2z^2)(\mathscr{A}_{22}y^2z^2)-(\mathscr{A}_{12}y^2z^2)(\mathscr{A}_{21}y^2z^2)>0$ can be checked by the minimum $M$-eigenvalue of the tensor $\mathscr{T}_A^1$.
\item It is easy to see that $\mbox{\rm det}(A(y,z))$ is a homogeneous polynomial of degree 12 with special structure. We can define the unique paired symmetric tensor $\mathscr{T}_A^2\in \mathbb{T}_{12,3}$ such that
    $$\mathscr{T}_A^2y^6z^6=\mbox{\rm det}(A(y,z)),$$
    where $(\mathscr{T}_A^2)_{i_1i_2\cdots i_{12}}=(\mathscr{T}_A^2)_{\sigma(i_1\cdots i_6)\sigma(i_7\cdots i_{12})}$ with $\sigma(i_1\cdots i_6)$ being an arbitrary permutation of $i_1i_2i_3i_4i_5i_6$. Thus, $\mbox{\rm det}(A(y,z))>0$ can be checked by the minimum $M$-eigenvalue of the tensor $\mathscr{T}_A^2$.
\end{itemize}
Similarly, we can define paired symmetric tensors $\mathscr{T}_B^1,\mathscr{T}_C^1\in \mathbb{T}_{8,3}$ and $\mathscr{T}_B^2,\mathscr{T}_C^2\in \mathbb{T}_{12,3}$ by
\begin{eqnarray*}
\begin{array}{l}
\mathscr{T}_B^1x^4z^4=(\mathscr{B}_{11}x^2z^2)(\mathscr{B}_{22}x^2z^2)-(\mathscr{B}_{12}x^2z^2)(\mathscr{B}_{21}x^2z^2),\quad \mathscr{T}_B^2x^6z^6=\mbox{\rm det}(B(x,z));\\
\mathscr{T}_C^1x^4y^4=(\mathscr{C}_{11}x^2y^2)(\mathscr{C}_{22}x^2y^2)-(\mathscr{C}_{12}x^2y^2)(\mathscr{C}_{21}x^2y^2),\quad \mathscr{T}_C^2x^6y^6=\mbox{\rm det}(C(x,y)).
\end{array}
\end{eqnarray*}
By combining Theorem \ref{thm-t-Meig-add} with Theorem \ref{thm-six-poly-1}, we have the following results.
\begin{Theorem}\label{thm-six-poly-2}
For any paired symmetric tensor $\mathscr{A}\in \mathbb{T}_{6,3}$, the polynomial $\mathscr{A}x^2y^2z^2$ defined by (\ref{e-6-poly}) is positive definite if and only if one of the following results holds.
\begin{itemize}
\item[(i)] The smallest $M$-eigenvalues of tensors $\mathscr{A}_{11}$, $\mathscr{T}_A^1$ and $\mathscr{T}_A^2$ are positive.
\item[(ii)] The smallest $M$-eigenvalues of tensors $\mathscr{B}_{11}$, $\mathscr{T}_B^1$ and $\mathscr{T}_B^2$ are positive.
\item[(iii)] The smallest $M$-eigenvalue of tensors $\mathscr{C}_{11}$, $\mathscr{T}_C^1$ and $\mathscr{T}_C^2$ are positive.
\end{itemize}
Furthermore, if $\mathscr{A}$ is an elasticity tensor, then the above (i), (ii) and (iii) are the same.
\end{Theorem}

\section{Sum of Squares and Positive Definiteness}

In this section, we investigate the SOS properties of polynomials defined by fourth order three dimensional paired symmetric (elasticity) tensors and sixth order three dimensional paired symmetric (elasticity) tensors, respectively. In particular, we give necessary and/or sufficient conditions of the concerned tensor being positive definite.

\noindent{\bf 5.1 Fourth order paired symmetric tensors}. For any paired symmetric tensor $\mathscr{A}=(a_{ijkl})\in \mathbb{T}_{4,3}$, we assume that the biquadratic form is defined by (\ref{e-4-poly}), i.e.,
$$
\mathscr{A}x^2y^2=\sum_{i,j,k,l=1}^3a_{ijkl}x_ix_jy_ky_l,\quad \forall x,y\in \mathbb{R}^3.
$$
We first investigate sufficient conditions under which the biquadratic form defined by (\ref{e-4-poly}) is an SOS of bilinear forms, or is positive definite.

\begin{Theorem}\label{thm-sos-suff}
Let $i_1i_2\cdots i_9$ and $j_1j_2\cdots j_9$ be two arbitrary permutations of $123123123$. For any paired symmetric (elasticity) tensor $\mathscr{A}\in \mathbb{T}_{4,3}$, let $M$ defined by Definition \ref{def-unfolded-m-4} be an unfolded matrix of $\mathscr{A}$ with respect to indices $i_1i_2\cdots i_9$ and $j_1j_2\cdots j_9$. Then, the following results hold.
\begin{itemize}
\item [(i)] If $M$ is  positive semidefinite, then the biquadratic form defined by (\ref{e-4-poly}) is an SOS of bilinear forms.
\item[(ii)] If $M$ is  positive definite, then the biquadratic form defined by (\ref{e-4-poly}) is positive definite.
\end{itemize}
\end{Theorem}

\noindent {\bf Proof}. (i) For any $x,y\in \mathbb{R}^3$, by given indices $i_1i_2\cdots i_9$ and $j_1j_2\cdots j_9$, we define $w\in \mathbb{R}^9$ by
\begin{eqnarray*}%\label{e-def-vector-u}
w:=(x_{i_1}y_{j_1},x_{i_2}y_{j_2},x_{i_3}y_{j_3},x_{i_4}y_{j_4},x_{i_5}y_{j_5},x_{i_6}y_{j_6},x_{i_7}y_{j_7},x_{i_8}y_{j_8},x_{i_9}y_{j_9})^\top.
\end{eqnarray*}
Then, from Definition \ref{def-unfolded-m-4}, it is not difficult to show that
\begin{eqnarray}\label{e-thm-sos-suff1}
\mathscr{A}x^2y^2=\sum_{i,j,k,l=1}^3a_{ijkl}x_ix_jy_ky_l=w^\top Mw, \quad \forall x,y\in \mathbb{R}^3,
\end{eqnarray}
and from Proposition \ref{pro-matrix-0} it follows that the matrix $M$ is symmetric.

Since the matrix $M$ is positive semidefinite, there are real numbers $\lambda_1,\lambda_2,\ldots,\lambda_9$ with $\lambda_1\geq\lambda_2\geq\ldots\geq\lambda_r\geq \lambda_{r+1}=\cdots=\lambda_9=0$  and an orthogonal matrix $Q=(q_1\; q_2\;\cdots\; q_9)$ such that
\begin{eqnarray}\label{e-thm-sos-suff2}
Q^\top MQ=\mbox{\rm diag}(\lambda_1, \lambda_2, \ldots, \lambda_9),\quad \mbox{\rm i.e.},\quad
M=\sum_{i=1}^r\lambda_iq_iq_i^\top.
\end{eqnarray}
This, together with (\ref{e-thm-sos-suff1}), implies that
\begin{eqnarray}\label{e-thm-sos-suff3}
\mathscr{A}x^2y^2=w^\top Mw=\sum_{i=1}^r\lambda_iw^\top q_iq_i^\top w=\sum_{i=1}^r\left(\sqrt{\lambda_i}\, w^\top q_i\right)^2,
\end{eqnarray}
which means that $\mathscr{A}x^2y^2$ is an SOS of bilinear forms.

(ii) Since the matrix $M$ is positive definite, there are real numbers $\lambda_1,\lambda_2,\ldots$, $\lambda_9$ with $\lambda_i>0$ for all $i\in \{1,2,\ldots,9\}$ and an orthogonal matrix $Q=(q_1\; q_2\;\cdots\; q_9)$ such that (\ref{e-thm-sos-suff2}) holds; and hence, (\ref{e-thm-sos-suff3}) holds. If $x\neq 0$ and $y\neq 0$, then it is easy to see that $w\neq 0$, which further implies that there exists at least an index $i\in \{1,2,\ldots,9\}$ such that $w^\top q_i\neq 0$ since the matrix $Q$ is invertible. Thus, when the matrix $M$ is  positive definite, it follows from (\ref{e-thm-sos-suff3}) that $\mathscr{A}x^2y^2>0$ for any $x\neq 0$ and $y\neq 0$, i.e., $\mathscr{A}x^2y^2$ is positive definite.
 \ep

Especially, we have following results.
\begin{Theorem}\label{thm-sos-suff-s}
For any paired symmetric tensor $\mathscr{A}\in \mathbb{T}_{4,3}$, let matrices $M^1$ and $M^2$ be defined by (\ref{e-matrix-A}) and (\ref{e-matrix-B}), respectively. Then, the following results hold.
\begin{itemize}
\item The biquadratic form defined by (\ref{e-4-poly}) is an SOS of bilinear forms if the matrix $M^1$ is  positive semidefinite or the matrix $M^2$ is positive semidefinite. Furthermore, if $\mathscr{A}\in \mathbb{T}_{4,3}$ is an elasticity tensor, then the results mentioned above are the same.
\item The biquadratic form defined by (\ref{e-4-poly}) is positive definite if the matrix $M^1$ is  positive definite or the matrix $M^2$ is positive definite. Furthermore, if $\mathscr{A}\in \mathbb{T}_{4,3}$ is an elasticity tensor, then the results mentioned above are the same.
\end{itemize}
\end{Theorem}

\noindent {\bf Proof}. For any $x,y\in \mathbb{R}^3$, we define two vectors $u,v\in \mathbb{R}^9$ by
\begin{eqnarray}\label{e-def-vector-u}
\begin{array}{rcl}
u&:=&(x_1y_1,x_1y_2,x_1y_3,x_2y_1,x_2y_2,x_2y_3,x_3y_1,x_3y_2,x_3y_3)^\top,\\
v&:=&(x_1y_1,x_2y_1,x_3y_1,x_1y_2,x_2y_2,x_3y_2,x_1y_3,x_2y_3,x_3y_3)^\top.
\end{array}
\end{eqnarray}
Then, it is easy to show that
$$
\mathscr{A}x^2y^2=\sum_{i,j,k,l=1}^3a_{ijkl}x_ix_jy_ky_l=u^\top M^1u=v^\top M^2v,\quad \forall x,y\in \mathbb{R}^3.
$$
Thus, by a similar way as Theorem \ref{thm-sos-suff} and using Proposition \ref{pro-matrix-1}, we can complete the proof. \ep

Next, we define a class of tensors with the help of $\mathscr{A}\in \mathbb{T}_{4,3}$.
\begin{Definition}\label{def-semi-p-s}
For any paired symmetric tensor $\mathscr{A}\in \mathbb{T}_{4,3}$, we define a tensor $\mathscr{B}\in \mathbb{T}_{4,3}$ whose entries satisfy
$$
b_{ijkl}=b_{jilk},\;\; b_{ijlk}=b_{jikl},\;\; \mbox{\rm and}\;\; b_{ijkl}+b_{jilk}+b_{ijlk}+b_{jikl}=4a_{ijkl}.
$$
We say that $\mathscr{B}$ is a {\bf semi-paired symmetric tensor} of $\mathscr{A}$.
\end{Definition}

Especially, we give the following remark.
\begin{Remark}\label{remark-def-semi-p-s}
We define two sets of indices by
\begin{eqnarray*}
\mathcal{I}&:=&\{(i,j,k,l)\in \{1,2,3\}^4: i\neq j,k\neq l\};\\
\mathcal{J}&:=&\{(i,j,k,l)\in \{1,2,3\}^4: i<j, l<k\;\;\mbox{\rm or}\;\; i>j, k>l\}.
\end{eqnarray*}
For any paired symmetric tensor $\mathscr{A}=(a_{ijkl})\in \mathbb{T}_{4,3}$, we define a tensor $\mathscr{B}=(b_{ijkl})\in \mathbb{T}_{4,3}$ whose entries satisfy
\begin{eqnarray*}%\label{e-def-semi-ps}
b_{ijkl}=\left\{\begin{array}{ll}
2a_{ijlk} & \mbox{\rm if}\; (i,j,l,k)\in \mathcal{J},\\
0         & \mbox{\rm if}\; (i,j,l,k)\in \mathcal{I}\setminus \mathcal{J},\\
a_{ijlk}  & \mbox{\rm otherwise},
\end{array}\right.  \quad \forall i,j,k,l\in \{1,2,3\},
\end{eqnarray*}
then $\mathscr{B}$ is a {\bf semi-paired symmetric tensor} of $\mathscr{A}$.
\end{Remark}

It is easy to see that the biquadratic forms defined by a paired symmetric tensor $\mathscr{A}\in \mathbb{T}_{4,3}$ and by its semi-paired symmetric tensor $\mathscr{B}$ are the same. Thus, the following results hold.

\begin{Theorem}\label{thm-sos-suff-nece}
For any paired symmetric (elasticity) tensor $\mathscr{A}\in \mathbb{T}_{4,3}$, let $\mathscr{B}$ defined by Definition \ref{def-semi-p-s} be a semi-paired symmetric tensor of $\mathscr{A}$. Then, the quadratic form $\mathscr{A}x^2y^2$ is an SOS of bilinear forms if and only if the quadratic form $\mathscr{B}x^2y^2$ is an SOS of bilinear forms; and tensor $\mathscr{A}$ is  positive (semidefinite) definite if and only if tensor $\mathscr{B}$ is positive (semidefinite) definite.
\end{Theorem}

Thus, SOS property and positive (semidefiniteness) definiteness of the paired symmetric tensor $\mathscr{A}\in \mathbb{T}_{4,3}$ can be studied by investigating SOS property and positive (semidefiniteness) definiteness of its semi-paired symmetric tensor.

By Theorems \ref{thm-sos-suff} and \ref{thm-sos-suff-nece}, we have the following results.

\begin{Theorem}\label{thm-sos-suff-add}
For any paired symmetric (elasticity) tensor $\mathscr{A}\in \mathbb{T}_{4,3}$, let $\mathscr{B}$ defined by Definition \ref{def-semi-p-s} be a semi-paired symmetric tensor of $\mathscr{A}$. Suppose that $i_1i_2\cdots i_9$ and $j_1j_2\cdots j_9$ are two arbitrary permutations of $123123123$; and $S$ is an unfolded matrix of tensor $\mathscr{B}$  with respect to indices $i_1i_2\cdots i_9$ and $j_1j_2\cdots j_9$. Then, the matrix $S$ is symmetric. Furthermore, if the matrix $S$ is  positive semidefinite, then the biquadratic form defined by (\ref{e-4-poly}) is an SOS of bilinear forms; and if the matrix $S$ is  positive definite, then the biquadratic form defined by (\ref{e-4-poly}) is positive definite.
\end{Theorem}

\begin{Corollary}\label{coro-sos-4}
For any paired symmetric (elasticity) tensor $\mathscr{A}\in \mathbb{T}_{4,3}$, suppose that $\mathscr{B}=(b_{ijkl})$ defined by Remark \ref{remark-def-semi-p-s} is a semi-paired symmetric tensor of $\mathscr{A}$, and the unfolded matrix of $\mathscr{B}$ is defined by
\begin{eqnarray}\label{e-matrix-A-0}
M:=(m_{st})\quad \mbox{\rm with}\; m_{st}=m_{3(i-1)+k,3(j-1)+l}=b_{ijkl}\; \mbox{\rm for any}\; i,j,k,l\in \{1,2,3\};
\end{eqnarray}
If the matrix $M$ is  positive semidefinite, then the biquadratic form defined by (\ref{e-4-poly}) is an SOS of bilinear forms; and if the matrix $M$ is  positive definite, then the biquadratic form defined by (\ref{e-4-poly}) is positive definite.
\end{Corollary}

In the following, we propose a necessary and sufficient condition under which a fourth order three dimensional paired symmetric (elasticity) tensor is positive semidefinite. For this purpose, we need the following lemma.

\begin{Lemma}\label{lem-sos-suff-nece-1}
For any paired symmetric (elasticity) tensor $\mathscr{A}\in \mathbb{T}_{4,3}$, suppose that the biquadratic form defined by (\ref{e-4-poly}) can be written as
\begin{eqnarray}\label{e-lem-sos-1}
\mathscr{A}x^2y^2=\sum_{i,j,k,l=1}^3a_{ijkl}x_ix_jy_ky_l=u^\top Mu,\quad \forall x,y\in \mathbb{R}^3,
\end{eqnarray}
where the vector $u$ is given by (\ref{e-def-vector-u}) and $M=(m_{st})\in \mathbb{R}^{9\times 9}$ is a symmetric matrix.
Then, the matrix $M$ is an unfolded matrix of some semi-paired symmetric tensor of $\mathscr{A}$.
\end{Lemma}

\noindent {\bf Proof}. We define a tensor $\mathscr{B}=(b_{ijkl})\in \mathbb{T}_{4,3}$ by
\begin{eqnarray}\label{e-lem-sos-2}
\left(\begin{array}{ccccccccc}
b_{1111} & b_{1112} & b_{1113} & b_{1211} & b_{1212} & b_{1213} & b_{1311} & b_{1312} & b_{1313} \\
b_{1121} & b_{1122} & b_{1123} & b_{1221} & b_{1222} & b_{1223} & b_{1321} & b_{1322} & b_{1323} \\
b_{1131} & b_{1132} & b_{1133} & b_{1231} & b_{1232} & b_{1233} & b_{1331} & b_{1332} & b_{1333} \\
b_{2111} & b_{2112} & b_{2113} & b_{2211} & b_{2212} & b_{2213} & b_{2311} & b_{2312} & b_{2313} \\
b_{2121} & b_{2122} & b_{2123} & b_{2221} & b_{2222} & b_{2223} & b_{2321} & b_{2322} & b_{2323} \\
b_{2131} & b_{2132} & b_{2133} & b_{2231} & b_{2232} & b_{2233} & b_{2331} & b_{2332} & b_{2333} \\
b_{3111} & b_{3112} & b_{3113} & b_{3211} & b_{3212} & b_{3213} & b_{3311} & b_{3312} & b_{3313} \\
b_{3121} & b_{3122} & b_{3123} & b_{3221} & b_{3222} & b_{3223} & b_{3321} & b_{3322} & b_{3323} \\
b_{3131} & b_{3132} & b_{3133} & b_{3231} & b_{3232} & b_{3233} & b_{3331} & b_{3332} & b_{3333}
\end{array}\right):=(m_{st}).
\end{eqnarray}
By using the symmetry of matrix $M$, it is easy to show that
\begin{eqnarray}\label{e-lem-sos-3}
b_{ijkl}=b_{jilk}\quad \mbox{\rm and}\quad b_{ijlk}=b_{jikl},\quad \forall i,j,k,l\in \{1,2,3\}.
\end{eqnarray}
Moreover, by combining (\ref{e-lem-sos-1}) with (\ref{e-lem-sos-2}) we have
$$
\mathscr{A}x^2y^2=\sum_{i,j,k,l=1}^3a_{ijkl}x_ix_jy_ky_l=u^\top Mu=\sum_{i,j,k,l=1}^3b_{ijkl}x_ix_jy_ky_l=\mathscr{B}x^2y^2,
$$
which, together with (\ref{e-lem-sos-3}) and the paired symmetry of $\mathscr{A}$, implies that
\begin{eqnarray}\label{e-lem-sos-4}
b_{ijkl}+b_{jilk}+b_{ijlk}+b_{jikl}=4a_{ijkl},\quad \forall i,j,k,l\in \{1,2,3\}.
\end{eqnarray}
Thus,  by Definition \ref{def-semi-p-s}, (\ref{e-lem-sos-3}) and (\ref{e-lem-sos-4}), it follows that $\mathscr{B}$ is a semi-paired symmetric tensor of $\mathscr{A}$.

In addition, from (\ref{e-lem-sos-2}) it is easy to see that
\begin{eqnarray*}
\begin{array}{l}
m_{st}=b_{i_si_tj_sj_t},\quad \forall s,t\in \{1,2,\ldots,9\}\;\;\mbox{\rm with}\\
\qquad (i_1,i_2,i_3,i_4,i_5,i_6,i_7,i_8,i_9)=(1,1,1,2,2,2,3,3,3),\\
\qquad (j_1,j_2,j_3,j_4,j_5,j_6,j_7,j_8,j_9)=(1,2,3,1,2,3,1,2,3).
\end{array}
\end{eqnarray*}
Thus, by Definition \ref{def-unfolded-m-4} we obtain that the matrix $M$ is an unfolded matrix of tensor $\mathscr{B}$ with respect to indices $111222333$ and $123123123$.
\ep

\begin{Theorem}\label{thm-sos-suff-nece-1}
For any paired symmetric (elasticity) tensor $\mathscr{A}\in \mathbb{T}_{4,3}$, the biquadratic form defined by (\ref{e-4-poly}) is an SOS of bilinear forms if and only if an unfolded matrix of some semi-paired symmetric tensor of $\mathscr{A}$ is  positive semidefinite.
\end{Theorem}

\noindent {\bf Proof}. The sufficiency can be obtained by Theorem \ref{thm-sos-suff-add}(i). We now show the necessity. Since $\mathscr{A}x^2y^2$ is an SOS of bilinear forms, there exists some positive integer $r\geq 1$ such that
$$
\mathscr{A}x^2y^2=\sum_{i,j,k,l=1}^3a_{ijkl}x_ix_jy_ky_l=\sum_{s=1}^r\left( \sum_{i,j=1}^3\alpha^s_{ij}x_iy_j\right)^2, \quad \forall x,y\in \mathbb{R}^3
$$
where $\alpha^s_{ij}\in \mathbb{R}$ for any $s\in \{1,2,\ldots,r\}$ and $i,j\in \{1,2,3\}$. For any $s\in \{1,2,\ldots,r\}$, we denote
$$
q_s:=(\alpha^s_{11},\alpha^s_{12},\alpha^s_{13},\alpha^s_{21},\alpha^s_{22},\alpha^s_{23},\alpha^s_{31},\alpha^s_{32},\alpha^s_{33})^\top.
$$
Let the vector $u\in \mathbb{R}^n$ be defined by (\ref{e-def-vector-u}). Then,
$$
\mathscr{A}x^2y^2=\sum_{s=1}^r\left(q_s^\top u\right)^2=u^\top\left(\sum_{s=1}^rq_sq_s^\top \right)u, \quad \forall x,y\in \mathbb{R}^3.
$$
By Lemma \ref{lem-sos-suff-nece-1}, we can obtain that the matrix $\sum_{s=1}^rq_sq_s^\top$ is an unfolded matrix of some semi-paired symmetric tensor of $\mathscr{A}$, which is positive semidefinite.
  \ep

In the following, we give a necessary condition of a fourth order three dimensional paired symmetric tensor being an SOS of bilinear forms, which is convenient for judging some fourth order three dimensional paired symmetric tensors not be an SOS.
\begin{Theorem}\label{thm-sos-necessary}
For any paired symmetric (elasticity) tensor $\mathscr{A}\in \mathbb{T}_{4,3}$, if the biquadratic form defined by (\ref{e-4-poly}) is an SOS of bilinear forms, then it follows that
$$
a_{iikk}\geq 0,\quad \forall i,k\in \{1,2,3\}.
$$
\end{Theorem}

\noindent {\bf Proof}. By the assumption, we can let
$$
\mathscr{A}x^2y^2=\sum_{i,j,k,l=1}^3a_{ijkl}x_ix_jy_ky_l=\sum_{s=1}^r\left( \sum_{i,k=1}^3\alpha^s_{ik}x_iy_k\right)^2, \quad \forall x,y\in \mathbb{R}^3
$$
where $r\geq 1$ is a positive integer and $\alpha^s_{ik}\in \mathbb{R}$ for any $s\in \{1,2,\ldots,r\}$ and $i,k\in \{1,2,3\}$. Then, it is easy to show that
\begin{eqnarray*}
a_{iikk}=\sum_{s=1}^r\left(\alpha^s_{ik}\right)^2,\quad \forall i,k\in \{1,2,3\},
\end{eqnarray*}
which completes the proof.    \ep

\noindent{\bf 5.2 Sixth order paired symmetric tensors}. In this part, we consider the homogeneous polynomial defined by a tensor $\mathscr{A}=(a_{ijklpq})\in \mathbb{T}_{6,3}$, which is given by (\ref{e-6-poly}).

Similar to Theorem \ref{thm-sos-suff}, we have the following results.
\begin{Theorem}\label{thm-sos-suff-6}
Let $i_1i_2\cdots i_{27}$ and $j_1j_2\cdots j_{27}$ be two arbitrary permutations of $\underbrace{123123\cdots123}_{27}$. For any paired symmetric (elasticity) tensor $\mathscr{A}\in \mathbb{T}_{6,3}$, suppose that $N$ defined by Definition \ref{def-unfolded-m-6} is an unfold matrix of $\mathscr{A}$ with respect to indices $i_1i_2\cdots i_{27}$ and $j_1j_2\cdots j_{27}$. Then, the following results hold.
\begin{itemize}
\item [(i)] If $N$ is  positive semidefinite, then the polynomial defined by (\ref{e-6-poly}) is an SOS of trilinear forms.
\item[(ii)] If $N$ is  positive definite, then the polynomial defined by (\ref{e-6-poly}) is positive definite.
\end{itemize}
\end{Theorem}

Especially, we have the following results.

\begin{Theorem}\label{thm-6-sos-suff-s}
For any paired symmetric tensor $\mathscr{A}=(a_{ijklpq})\in \mathbb{T}_{6,3}$, let matrices $N^1, N^2, \ldots, N^6$ be defined by (\ref{e-6-matrix-1}). Then, the following results hold.
\begin{itemize}
\item[(i)] The polynomial defined by (\ref{e-6-poly}) is an SOS of trilinear polynomials if one of matrices $N^1, N^2, \ldots, N^6$ is positive semidefinite. Furthermore, if $\mathscr{A}=(a_{ijklpq})\in \mathbb{T}_{6,3}$ is an elasticity tensor, then the results mentioned above are the same.
\item[(ii)] The polynomial defined by (\ref{e-6-poly}) is positive definite if one of matrices $N^1$, $N^2$, $\ldots$, $N^6$ is positive definite. Furthermore, if $\mathscr{A}=(a_{ijklpq})\in \mathbb{T}_{6,3}$ is an elasticity tensor, then the results mentioned above are the same.
\end{itemize}
\end{Theorem}

\noindent {\bf Proof}. For any $x,y,z\in \mathbb{R}^3$, we define six vectors $u^i\in \mathbb{R}^{27}\; (i\in \{1,2,\ldots,6\})$ by
\begin{eqnarray*}
\begin{array}{rcl}
(u^1)^\top&:=&(x_1y_1z_1,x_1y_1z_2,x_1y_1z_3,x_1y_2z_1,x_1y_2z_2,x_1y_2z_3,x_1y_3z_1,x_1y_3z_2,x_1y_3z_3,\\
      &  &\; x_2y_1z_1,x_2y_1z_2,x_2y_1z_3,x_2y_2z_1,x_2y_2z_2,x_2y_2z_3,x_2y_3z_1,x_2y_3z_2,x_2y_3z_3,\\
      &  &\; x_3y_1z_1,x_3y_1z_2,x_3y_1z_3,x_3y_2z_1,x_3y_2z_2,x_3y_2z_3,x_3y_3z_1,x_3y_3z_2,x_3y_3z_3);\\
(u^2)^\top&:=&(x_1y_1z_1,x_1y_1z_2,x_1y_1z_3,x_2y_1z_1,x_2y_1z_2,x_2y_1z_3,x_3y_1z_1,x_3y_1z_2,x_3y_1z_3,\\
      &  &\; x_1y_2z_1,x_1y_2z_2,x_1y_2z_3,x_2y_2z_1,x_2y_2z_2,x_2y_2z_3,x_3y_2z_1,x_3y_2z_2,x_3y_2z_3,\\
      &  &\; x_3y_1z_1,x_3y_1z_2,x_3y_1z_3,x_3y_2z_1,x_3y_2z_2,x_3y_2z_3,x_3y_3z_1,x_3y_3z_2,x_3y_3z_3);\\
(u^3)^\top&:=&(x_1y_1z_1,x_2y_1z_1,x_3y_1z_1,x_1y_2z_1,x_2y_2z_1,x_3y_2z_1,x_1y_3z_1,x_2y_3z_1,x_3y_3z_1,\\
      &  &\; x_1y_1z_2,x_2y_1z_2,x_3y_1z_2,x_1y_2z_2,x_2y_2z_2,x_3y_2z_2,x_1y_3z_2,x_2y_3z_2,x_3y_3z_2,\\
      &  &\; x_1y_1z_3,x_2y_1z_3,x_3y_1z_3,x_1y_2z_3,x_2y_2z_3,x_3y_2z_3,x_1y_3z_3,x_2y_3z_3,x_3y_3z_3);\\
(u^4)^\top&:=&(x_1y_1z_1,x_1y_2z_1,x_1y_3z_1,x_1y_1z_2,x_1y_2z_2,x_1y_3z_2,x_1y_1z_3,x_1y_2z_3,x_1y_3z_3,\\
      &  &\; x_2y_1z_1,x_2y_2z_1,x_2y_3z_1,x_2y_1z_2,x_2y_2z_2,x_2y_3z_2,x_2y_1z_3,x_2y_2z_3,x_2y_3z_3,\\
      &  &\; x_3y_1z_1,x_3y_1z_2,x_3y_1z_3,x_3y_2z_1,x_3y_2z_2,x_3y_2z_3,x_3y_3z_1,x_3y_3z_2,x_3y_3z_3);\\
(u^5)^\top&:=&(x_1y_1z_1,x_1y_2z_1,x_1y_3z_1,x_2y_1z_1,x_2y_2z_1,x_2y_3z_1,x_3y_1z_1,x_3y_2z_1,x_3y_3z_1,\\
      &  &\; x_1y_1z_2,x_1y_2z_2,x_1y_3z_2,x_2y_1z_2,x_2y_2z_2,x_2y_3z_2,x_3y_1z_2,x_3y_2z_2,x_3y_3z_2,\\
      &  &\; x_1y_1z_3,x_2y_1z_3,x_3y_1z_3,x_1y_2z_3,x_2y_2z_3,x_3y_2z_3,x_1y_3z_3,x_2y_3z_3,x_3y_3z_3);\\
(u^6)^\top&:=&(x_1y_1z_1,x_2y_1z_1,x_3y_1z_1,x_1y_1z_2,x_2y_1z_2,x_3y_1z_2,x_1y_1z_3,x_2y_1z_3,x_3y_1z_3,\\
      &  &\; x_1y_2z_1,x_2y_2z_1,x_3y_2z_1,x_1y_2z_2,x_2y_2z_2,x_3y_2z_2,x_1y_2z_3,x_2y_2z_3,x_3y_2z_3,\\
      &  &\; x_1y_3z_1,x_1y_3z_2,x_1y_3z_3,x_2y_3z_1,x_2y_3z_2,x_2y_3z_3,x_3y_3z_1,x_3y_3z_2,x_3y_3z_3).
\end{array}
\end{eqnarray*}
With the help of the definitions of matrices $N^1, N^2, \ldots, N^6$ by (\ref{e-6-matrix-1}), it is easy to verify that
\begin{eqnarray*}%\label{e-thm-sos-1}
\mathscr{A}x^2y^2z^2=\sum_{i,j,k,l,p,q=1}^3a_{ijklpq}x_ix_jy_ky_lz_pz_q
=(u^1)^\top N^1u^1=(u^2)^\top N^2u^2=\cdots=(u^6)^\top N^6u^6.
\end{eqnarray*}
Thus, similar to Theorem \ref{thm-sos-suff-s}, we can obtain the desired results.  \ep

Similar to Definition \ref{def-semi-p-s}, we give the concept of semi-paired symmetric tensor of $\mathscr{A}\in \mathbb{T}_{6,3}$.
\begin{Definition}\label{def-semi-p-s-6}
For any paired symmetric tensor $\mathscr{A}\in \mathbb{T}_{6,3}$, we define a tensor $\mathscr{B}\in \mathbb{T}_{6,3}$ whose entries satisfy
\begin{eqnarray*}
\begin{array}{l}
b_{ijklpq}=b_{jilkqp},\;\; b_{jiklpq}=b_{ijlkqp},\;\; b_{ijlkpq}=b_{jiklqp},\;\; b_{ijklqp}=b_{jilkpq},\\
\mbox{\rm and}\;\; b_{ijklpq} + b_{jiklpq}+b_{ijlkpq}+b_{ijklqp}=4a_{ijklpq}.
\end{array}
\end{eqnarray*}
We say that $\mathscr{B}$ is a {\bf semi-paired symmetric tensor} of $\mathscr{A}$.
\end{Definition}

Especially, we give the following remark.
\begin{Remark}\label{remark-def-semi-p-s-6}
We define several index sets by
\begin{eqnarray*}%\label{e-index-partition-0}
\begin{array}{rcl}
\mathcal{T}_1&:=&\left\{(i,j,k,l,p,q)\in \{1,2,3\}^6: i\neq j, k\neq l\; \mbox{\rm and}\; p\neq q\right\},\vspace{2mm}\\
\mathcal{T}_2&:=&\left\{(i,j,k,l,p,q)\in \{1,2,3\}^6: \begin{array}{l}
i\neq j, k\neq l\; \mbox{\rm and}\; p=q; \\
\mbox{\rm or}\; i\neq j, k=l\; \mbox{\rm and}\; p\neq q; \\
\mbox{\rm or}\; i=j, k\neq l\; \mbox{\rm and}\; p\neq q\end{array} \right\},\vspace{2mm}\\
\mathcal{T}_3&:=&\left\{(i,j,k,l,p,q)\in \{1,2,3\}^6: \begin{array}{l}
i\neq j, k=l\; \mbox{\rm and}\; p=q; \\
\mbox{\rm or}\; i=j, k\neq l\; \mbox{\rm and}\; p=q; \\
\mbox{\rm or}\; i=j, k=l\; \mbox{\rm and}\; p\neq q\end{array} \right\},\vspace{2mm}\\
\mathcal{T}_4&:=& \left\{(i,j,k,l,p,q)\in \{1,2,3\}^6: i=j, k=l\; \mbox{\rm and}\; p=q\right\}
\end{array}
\end{eqnarray*}
and
\begin{eqnarray*}%\label{e-index-partition}
\begin{array}{rcl}
\mathcal{S}_1&:=&\left\{(i,j,k,l,p,q)\in \mathcal{T}_1: \begin{array}{l}
i<j, k<l\; \mbox{\rm and}\; p<q;\\
\mbox{\rm or}\; i>j, k>l\; \mbox{\rm and}\; p>q \end{array} \right\},\vspace{2mm}\\
\mathcal{S}_2&:=&\left\{(i,j,k,l,p,q)\in \mathcal{T}_2: \begin{array}{l}
i<j, k<l\; \mbox{\rm and}\; p=q; \\
\mbox{\rm or}\; i>j, k>l\; \mbox{\rm and}\; p=q \end{array} \right\},\vspace{2mm}\\
\mathcal{S}_3&:=&\left\{(i,j,k,l,p,q)\in \mathcal{T}_2: \begin{array}{l}
i<j, k=l\; \mbox{\rm and}\; p<q; \\
\mbox{\rm or}\; i>j, k=l\; \mbox{\rm and}\; p>q\end{array} \right\},\vspace{2mm}\\
\mathcal{S}_4&:=&\left\{(i,j,k,l,p,q)\in \mathcal{T}_2: \begin{array}{l}
i=j, k<l\; \mbox{\rm and}\; p<q; \\
\mbox{\rm or}\; i=j, k>l\; \mbox{\rm and}\; p>q\end{array} \right\}.\vspace{2mm}\\
\end{array}
\end{eqnarray*}
For any paired symmetric tensor $\mathscr{A}\in \mathbb{T}_{6,3}$, we define a tensor $\mathscr{B}=(b_{ijklpq})\in \mathbb{T}_{6,3}$ whose entries satisfy
\begin{eqnarray*}%\label{e-def-semi-ps}
b_{ijklpq}=\left\{\begin{array}{ll}
4a_{ijlkpq} & \mbox{\rm if}\; (i,j,l,k,p,q)\in \mathcal{S}_1,\\
0           & \mbox{\rm if}\; (i,j,l,k,p,q)\in \mathcal{T}_1\setminus \mathcal{S}_1,\\
2a_{ijlkpq} & \mbox{\rm if}\; (i,j,l,k,p,q)\in \mathcal{S}_2\bigcup \mathcal{S}_3\bigcup \mathcal{S}_4,\\
0           & \mbox{\rm if}\; (i,j,l,k,p,q)\in \mathcal{T}_2\setminus \{\mathcal{S}_2\bigcup \mathcal{S}_3\bigcup \mathcal{S}_4\},\\
a_{ijlkpq}  & \mbox{\rm otherwise},
\end{array}\right.  \quad \forall i,j,k,l,p,q\in \{1,2,3\},
\end{eqnarray*}
then $\mathscr{B}$ is a {\bf semi-paired symmetric tensor} of $\mathscr{A}$.
\end{Remark}

Similar to Theorem \ref{thm-sos-suff-nece}, we have the following results.
\begin{Theorem}\label{thm-6-sos-suff-nece}
For any paired symmetric (elasticity) tensor $\mathscr{A}\in \mathbb{T}_{6,3}$, suppose that $\mathscr{B}$ defined by Definition \ref{def-semi-p-s-6} is an semi-paired symmetric tensor of $\mathscr{A}$. Then, the polynomial $\mathscr{A}x^2y^2z^2$ is an SOS of trilinear forms if and only if the polynomial $\mathscr{B}x^2y^2z^2$ is an SOS of trilinear forms; and tensor $\mathscr{A}$ is  positive (semidefinite) definite if and only if tensor $\mathscr{B}$ is positive (semidefinite) definite.
\end{Theorem}

\begin{Corollary}\label{coro-sos-6}
For any paired symmetric (elasticity) tensor $\mathscr{A}\in \mathbb{T}_{6,3}$, suppose that $\mathscr{B}=(b_{ijkl})$ defined by Remark \ref{remark-def-semi-p-s-6} is a semi-paired symmetric tensor of $\mathscr{A}$, and the unfolded matrix of $\mathscr{B}$ is defined by
\begin{eqnarray}\label{e-matrix-N-0}
\begin{array}{l}
N:=(n_{st})\;\; \mbox{\rm with}\\
\quad  n_{st}=n_{3[3(i-1)+(k-1)]+p,3[3(j-1)+(l-1)]+q}=b_{ijklpq},\; \forall i,j,k,l,p,q\in \{1,2,3\};
\end{array}
\end{eqnarray}
If the matrix $N$ is  positive semidefinite, then the polynomial defined by (\ref{e-6-poly}) is an SOS of trilinear forms; and if the matrix $N$ is  positive definite, then the polynomial defined by (\ref{e-6-poly}) is positive definite.
\end{Corollary}

Similar to Theorem \ref{thm-sos-suff-nece-1}, we have the following results.
\begin{Theorem}\label{thm-6-sos-suff-nece-1}
For any paired symmetric (elasticity) tensor $\mathscr{A}\in \mathbb{T}_{6,3}$, the polynomial defined by (\ref{e-6-poly}) is an SOS of trilinear forms if and only if an unfolded matrix of some semi-paired symmetric tensor of $\mathscr{A}$ is  positive semidefinite.
\end{Theorem}

The following necessary condition can be obtained in a similar way as Theorem \ref{thm-sos-necessary}.
\begin{Theorem}\label{thm-6-sos-necessary}
For any paired symmetric (elasticity) tensor $\mathscr{A}=(a_{ijklpq})\in \mathbb{T}_{6,3}$. If the homogeneous polynomial defined by (\ref{e-6-poly}) is an SOS of trilinear polynomials, then it follows that
$$
a_{iikkpp}\geq 0,\quad \forall i,k,p\in \{1,2,3\}.
$$
\end{Theorem}

\section{Algorithm for the Smallest $M$-Eigenvalue}

In order to judge positive definiteness of a fourth order three dimensional paired symmetric (elasticity) tensor, from Corollary \ref{thm-4-positive-basic} we may compute the smallest $M$-eigenvalue of the concerned tensor. Some methods can be applied to do it, such as the one given in \cite{NZ-16}. In this section, by using the special structure of the paired symmetric (elasticity) tensor, we propose a sequential sedefinite programming method for computing the smallest $M$-eigenvalue of a fourth order three dimensional paired symmetric (elasticity) tensor, by which we can judge whether a fourth order three dimensional paired symmetric (elasticity) tensor is positive definite or not.

Let $\mathscr{A}x^2y^2$ be defined by (\ref{e-4-poly}) and $r$ be some nonnegative integer, we define a polynomial $F_{r,s}: \mathbb{R}^3\times \mathbb{R}^3\rightarrow \mathbb{R}$ by
\begin{eqnarray}\label{e-4-poly-r}
F_{r,s}(x,y):=\left(\sum_{i=1}^3x_i^2\right)^r\left(\sum_{i=1}^3y_i^2\right)^s\mathscr{A}x^2y^2,
\end{eqnarray}
which is a homogeneous polynomial with deg$(F_r)=2(r+s)+4$. The following result can be obtained in a similar way as those in \cite[Corollary]{Reznick-95}, \cite[Theorem 3.2]{HHQ-13} and \cite[Lemma 3.5]{LNQY-09}.

\begin{Theorem}\label{thm-4-poly-r-sos1}
For any paired symmetric tensor $\mathscr{A}=(a_{ijkl})\in \mathbb{T}_{4,3}$, let $\mbox{\rm int}\mathcal{C}$ and $F_{r,s}$ be defined by (\ref{e-4-cone-1}) and (\ref{e-4-poly-r}), respectively. If $\mathscr{A}\in \mbox{\rm int}\mathcal{C}$, then it follows that $F_{r,s}(x,y)$ is an SOS for some sufficiently large integers $r,s\geq 0$.
\end{Theorem}

With Theorem \ref{thm-4-poly-r-sos1}, we define
\begin{eqnarray}\label{e-4-poly-r-cone}
\mathcal{K}:=\{\mathscr{A}\in \mathbb{T}_{4,3}: F_{r,s}(x,y)\;\mbox{\rm is an SOS for some}\; r,s\geq 0\}.
\end{eqnarray}
Then, similar to Theorem 3.3 in \cite{HHQ-13}, we can obtain the following result.

\begin{Theorem}\label{thm-4-poly-r-sos2}
For any paired symmetric tensor $\mathscr{A}=(a_{ijkl})\in \mathbb{T}_{4,3}$, $\mathcal{K}$ defined by (\ref{e-4-poly-r-cone}) is a convex cone. In particular, it follows that $\mbox{\rm int}\mathcal{C}\subseteq \mathcal{K}\subseteq \mathcal{C}$, where $\mathcal{C}$ and $\mbox{\rm int}\mathcal{C}$ are defined by (\ref{e-4-cone}) and (\ref{e-4-cone-1}), respectively.
\end{Theorem}

For any nonnegative integer $r$, by the multinomial theorem we have
$$
\left(\sum_{i=1}^3x_i^2\right)^r
=\sum_{r_1+r_2+r_3=r}\frac{r!}{r_1!\times r_2!\times r_3!}x_1^{2r_1}x_2^{2r_2}x_3^{2r_3}.
$$
Thus, the polynomial $F_{r,s}(x,y)$ defined by (\ref{e-4-poly-r}) is a linear combination of monomials
$$
x_1^{2r_1}x_2^{2r_2}x_3^{2r_3}y_1^{2s_1}y_2^{2s_2}y_3^{2s_3}x_ix_jy_ky_l
$$
for $i,j,k,l\in \{1,2,3\}$ and nonnegative integers $r_1,r_2,r_3,s_1,s_2,s_3$ satisfying $\sum_{i=1}^3r_i=r$ and $\sum_{i=1}^3s_i=s$. We assume that the corresponding vector of monomials is denoted by $v_{r+s+2}(x,y)$, which is first ordered in descending order of powers of $x_1,x_2,x_3$ and is then ordered in descending order of powers of $y_1,y_2,y _3$ given by
\begin{eqnarray}\label{e-4-poly-r-monomial}
\begin{array}{rcl}
v_{r+s+2}(x,y)^\top:&=&(x_1^{r+1}y_1^{s+1},x_1^{r+1}y_1^sy_2,\ldots, x_1^{r+1}y_3^{s+1}, x_1^rx_2y_2^{r+1}, \ldots, x_3^{r+1}y_3^{s+1}).
\end{array}
\end{eqnarray}
Then, the following result is immediate.

\begin{Theorem}\label{thm-4-poly-r-sos3}
For any paired symmetric tensor $\mathscr{A}=(a_{ijkl})\in \mathbb{T}_{4,3}$ and nonnegative integers $r$ and $s$, suppose that the homogeneous polynomial $F_{r,s}$ and the vector of monomials $v_{r+s+2}(x,y)$ are defined by (\ref{e-4-poly-r}) and (\ref{e-4-poly-r-monomial}), respectively. Then, the polynomial $F_{r,s}(x,y)$ is an SOS if and only if $F_{r,s}(x,y)=v_{r+s+2}(x,y)^\top H v_{r+s+2}(x,y)$, where $H$ is a positive semidefinite matrix.
\end{Theorem}

For any given nonnegative integers $r$ and $s$, we let $u_{r+s+2}(x,y)$ denote the vector made up of different monomials in $v_{r+s+2}(x,y)^\top v_{r+s+2}(x,y)$, which is first ordered in descending order of powers of $x_1,x_2,x_3$ and then ordered in descending order of powers of $y_1,y_2,y _3$; and denote the dimensions of vectors $v_{r+s+2}(x,y)$ and $u_{r+s+2}(x,y)$ by $d_{vrs}:=\mbox{\rm dim}(v_{r+s+2}(x,y))$ and $d_{urs}:=\mbox{\rm dim}(u_{r+s+2}(x,y))$, respectively.
\begin{itemize}
\item For $r=0$ and $s=0$, the vector $v_{r+s+2}(x,y)$ contains monomials
$$
x_iy_j,\quad \forall i,j\in \{1,2,3\},
$$
and hence, $d_{v00}=C_3^1C_3^1=9$; and the vector $u_{r+s+2}(x,y)$ contains monomials
$$
x_i^2y_ky_l\,(k<l),\; x_i^2y_k^2,\; x_ix_jy_ky_l\,(i<j,k<l),\; x_ix_jy_k^2\,(i<j),\quad \forall i,j,k,l\in \{1,2,3\},
$$
and hence, $d_{u00}=C_3^1C_3^2+C_3^1C_3^1+C_3^2C_3^2+C_3^1C_3^2=36$.
\item For $r=1$ and $s=0$, the vector $v_{r+s+2}(x,y)$ contains monomials
$$
x_1^{\alpha_1}x_2^{\alpha_2}x_3^{\alpha_3}y_1^{\beta_1}y_2^{\beta_2}y_3^{\beta_3},
$$
where $\alpha_i,\beta_i\geq0$ for all $i\in \{1,2,3\}$, $\alpha_1+\alpha_2+\alpha_3=2$ and $\beta_1+\beta_2+\beta_3=1$,
and hence, $d_{v10}=C_{3+2-1}^2C_{3+1-1}^1=18$; and the vector $u_{r+s+2}(x,y)$ contains monomials
$$
x_1^{\alpha_1}x_2^{\alpha_2}x_3^{\alpha_3}y_1^{\beta_1}y_2^{\beta_2}y_3^{\beta_3},
$$
where $\alpha_i,\beta_i\geq0$ for all $i\in \{1,2,3\}$, $\alpha_1+\alpha_2+\alpha_3=4$ and $\beta_1+\beta_2+\beta_3=2$, and hence,  $d_{u10}=C_{3+4-1}^4C_{3+2-1}^2=90$.
\item For $r=0$ and $s=1$, the vector $v_{r+s+2}(x,y)$ contains monomials
$$
x_1^{\alpha_1}x_2^{\alpha_2}x_3^{\alpha_3}y_1^{\beta_1}y_2^{\beta_2}y_3^{\beta_3},
$$
where $\alpha_i,\beta_i\geq0$ for all $i\in \{1,2,3\}$, $\alpha_1+\alpha_2+\alpha_3=1$ and $\beta_1+\beta_2+\beta_3=2$,
and hence, $d_{v01}=C_{3+1-1}^1C_{3+2-1}^2=18$; and the vector $u_{r+s+2}(x,y)$ contains monomials
$$
x_1^{\alpha_1}x_2^{\alpha_2}x_3^{\alpha_3}y_1^{\beta_1}y_2^{\beta_2}y_3^{\beta_3},
$$
where $\alpha_i,\beta_i\geq0$ for all $i\in \{1,2,3\}$, $\alpha_1+\alpha_2+\alpha_3=2$ and $\beta_1+\beta_2+\beta_3=4$, and hence,  $d_{u01}=C_{3+2-1}^2C_{3+4-1}^4=90$.
\item For $r=1$ and $s=1$, the vector $v_{r+s+2}(x,y)$ contains monomials
$$
x_1^{\alpha_1}x_2^{\alpha_2}x_3^{\alpha_3}y_1^{\beta_1}y_2^{\beta_2}y_3^{\beta_3},
$$
where $\alpha_i,\beta_i\geq0$ for all $i\in \{1,2,3\}$, $\alpha_1+\alpha_2+\alpha_3=2$ and $\beta_1+\beta_2+\beta_3=2$,
and hence, $d_{v01}=C_{3+2-1}^2C_{3+2-1}^2=36$; and the vector $u_{r+s+2}(x,y)$ contains monomials
$$
x_1^{\alpha_1}x_2^{\alpha_2}x_3^{\alpha_3}y_1^{\beta_1}y_2^{\beta_2}y_3^{\beta_3},
$$
where $\alpha_i,\beta_i\geq0$ for all $i\in \{1,2,3\}$, $\alpha_1+\alpha_2+\alpha_3=4$ and $\beta_1+\beta_2+\beta_3=4$, and hence,  $d_{u01}=C_{3+4-1}^4C_{3+4-1}^4=225$.
\end{itemize}
Furthermore, we introduce two operators $\mathcal{V}$ and $\mathcal{W}$ in the following:
\begin{itemize}
\item[$\mathcal{V}$:] For any nonnegative integers $r$ and $s$, we define an operator $\mathcal{V}: \mathbb{R}^{d_{vrs}\times d_{vrs}}\rightarrow \mathbb{R}^{d_{urs}}$ such that $[\mathcal{V}(H)]_i$ is the coefficient of the $i$th monomial in the vector $u_{r+s+2}(x,y)$ of the polynomial $v_{r+s+2}(x,y)^\top Hv_{r+s+2}(x,y)$ for any $H\in \mathbb{R}^{d_{vrs}\times d_{vrs}}$.
\item[$\mathcal{W}$:] For any nonnegative integers $r$ and $s$, we define an operator $\mathcal{W}$ which maps a homogeneous polynomial $F_{r,s}(x,y)$ defined by (\ref{e-4-poly-r}) to a vector in $\mathbb{R}^{d_{urs}}$ satisfying $F_{r,s}(x,y)=u_{r+s+2}(x,y)^\top \mathcal{W}(F_{r,s}(x,y))$.
\end{itemize}
Thus, by Theorem \ref{thm-4-poly-r-sos2} it follows that the set $\mathcal{K}$ defined by (\ref{e-4-poly-r-cone}) can be written as
\begin{eqnarray}\label{e-4-poly-r-cone-equi}
\mathcal{K}:=\left\{\mathscr{A}\in \mathbb{T}_{4,3}: \mathcal{W}(F_{r,s}(x,y))=\mathcal{V}(Q), r,s\in \{0,1,2,\ldots\}, Q\in \mathcal{S}_+^{d_{vrs}}\right\}.
\end{eqnarray}
For any $r,s\in \{0,1,2,\ldots\}$, we define
\begin{eqnarray}\label{e-4-poly-r-cone-equi-1}
\mathcal{K}_{r+s}:=\left\{\mathscr{A}\in \mathbb{T}_{4,3}: \mathcal{W}(F_{r,s}(x,y))=\mathcal{V}(Q), Q\in \mathcal{S}_+^{d_{vrs}}\right\}.
\end{eqnarray}
Then, the following result is immediate.

\begin{Theorem}\label{thm-4-poly-r-sos4}
For any paired symmetric tensor $\mathscr{A}=(a_{ijkl})\in \mathbb{T}_{4,3}$ and nonnegative integers $r$ and $s$, suppose that the set $\mathcal{K}_{r+s}$ is defined by (\ref{e-4-poly-r-cone-equi-1}). Then,
$$
\mathcal{K}_{r+s}\subseteq K_{r+s+1}\; \mbox{\rm for any}\; r,s\in \{0,1,2,\ldots\}\quad \mbox{\rm and}\quad \lim\limits_{r+s\rightarrow\infty} \mathcal{K}_{r+s}=\cup_{r+s=0}^\infty \mathcal{K}_{r+s}=\mathcal{K}.
$$
\end{Theorem}

For any paired symmetric tensor $\mathscr{A}=(a_{ijkl})\in \mathbb{T}_{4,3}$, it follows from Definition \ref{def-t-M-eigen} and Theorem \ref{thm-t-Meig-exist} that the smallest $M$-eigenvalue of $\mathscr{A}$ is the optimal value of the following biquadratic optimization over unit spheres \cite{LNQY-09}:
%\begin{eqnarray}\label{e-M-opt}
%\min \mathscr{A}x^2y^2\quad \mbox{\rm s.t.} \;\; x^\top x=1,\; y^\top y=1,
%\end{eqnarray}
\begin{eqnarray}\label{e-M-opt}
\begin{array}{cl}
\min & \mathscr{A}x^2y^2\\
\mbox{\rm s.t.} & x^\top x=1,\; y^\top y=1,
\end{array}
\end{eqnarray}
which is equivalent to
%\begin{eqnarray}\label{e-M-opt-1}
%\max \gamma\quad \mbox{\rm s.t.} \;\; \mathscr{A}x^2y^2\geq \gamma, \forall (x,y)\in \{(x,y)\in \mathbb{R}^3\times \mathbb{R}^3: x^\top x=1,\; y^\top y=1\}.
%\end{eqnarray}
\begin{eqnarray}\label{e-M-opt-1}
\begin{array}{cl}
\max & \gamma\\
\mbox{\rm s.t.} & \mathscr{A}x^2y^2\geq \gamma, \forall (x,y)\in \{(x,y)\in \mathbb{R}^3\times \mathbb{R}^3: x^\top x=1,\; y^\top y=1\}.
\end{array}
\end{eqnarray}
Let tensor $\mathscr{E}=(e_{ijkl})\in \mathbb{T}_{4,3}$ be defined by (\ref{e-partial-identity}).
Then, it is easy to see that problem (\ref{e-M-opt-1}) is equivalent to
%\begin{eqnarray}\label{e-M-opt-2}
%\max \gamma\quad \mbox{\rm s.t.} \;\; \mathscr{A}+\gamma\mathscr{E}\in \mathcal{K}.
%\end{eqnarray}
\begin{eqnarray}\label{e-M-opt-2}
\begin{array}{cl}
\max & \gamma\\
\mbox{\rm s.t.} & \mathscr{A}+\gamma\mathscr{E}\in \mathcal{K}.
\end{array}
\end{eqnarray}
Furthermore, by replacing the constraint $\mathscr{A}+\gamma\mathscr{E}\in \mathcal{K}$ by $\mathscr{A}+\gamma\mathscr{E}\in \mathcal{K}_{r+s}$, we obtain the $r+s$th order relaxation problem of (\ref{e-M-opt-2}), which can be written as
%\begin{eqnarray}\label{e-M-opt-3}
%\max \gamma\quad
%\mbox{\rm s.t.} \;\; \mathcal{V}(Q)=\mathcal{W}(F_{r,s}(x,y))+\gamma\mathcal{W}\left((\sum_{i=1}^3x_i^2)^{r+1}(\sum_{i=1}^3y_i^2)^{s+1}\right),\;
%                 Q\in \mathcal{S}_+^{d_{vrs}},
%\end{eqnarray}
\begin{eqnarray}\label{e-M-opt-3}
\begin{array}{cl}
\max & \gamma\\
\mbox{\rm s.t.} & \mathcal{V}(Q)=\mathcal{W}(F_{r,s}(x,y))+\gamma\mathcal{W}\left((\sum_{i=1}^3x_i^2)^{r+1}(\sum_{i=1}^3y_i^2)^{s+1}\right),\\
                & Q\in \mathcal{S}_+^{d_{vrs}},
\end{array}
\end{eqnarray}
where $r,s\in \{0,1,2,\ldots\}$. For any $r,s\in \{0,1,2,\ldots\}$, (\ref{e-M-opt-3}) is a semidefinite programming problem, which is denoted by $SDP(r+s)$. Thus, we can solve a sequence of $SDP(r+s)$ through increasing $r+s$ to obtain an approximation optimal solution of (\ref{e-M-opt}) up to a priori precision.

\section{Numerical Results}

In this section, we present preliminary numerical results to judge whether a fourth order three dimensional or a sixth order three dimensional paired symmetric tensor is positive definite or not. All numerical experiments were run in Matlab on a PC with 2.93 GHz CPU and 2.00 GB of RAM. We divide our experiments into the following two parts.

\noindent {\bf Part 1}. In this part, by using the sequential semidefinite programming method proposed in Section 6, we present preliminary numerical results for computing the smallest $M$-eigenvalue of the fourth order three dimensional paired symmetric tensor.  We use SDPT3 \cite{TTT-99} to solve the resulted conic linear programming problem.

\begin{Example}\label{exam1}
Let $\mathscr{A}=(a_{ijkl})\in \mathbb{T}_{4,3}$ be given by
\begin{eqnarray}\label{e-exam1-t-1}
\begin{array}{l}
a_{1111}=1,\; a_{2222}=1,\; a_{3333}=1,\; a_{1122}=2,\; a_{2233}=2,\; a_{3311}=2,\\
a_{1212}=a_{1221}=a_{2112}=a_{2121}=-\frac 12,\; a_{2323}=a_{2332}=a_{3223}=a_{3232}=-\frac 12,\\
a_{3131}=a_{3113}=a_{1331}=a_{1313}=-\frac 12.
\end{array}
\end{eqnarray}
\end{Example}

The biquadratic form defined by this tensor is
\begin{eqnarray}\label{e-exam1-p-2}
\begin{array}{rcl}
f(x,y)&=&x_1^2y_1^2+x_2^2y_2^2+x_3^2y_3^2-2(x_1x_2y_1y_2+x_2x_3y_2y_3+x_3x_1y_3y_1)\\
& &+2(x_1^2y_2^2+x_2^2y_3^2+x_3^2y_1^2),
\end{array}
\end{eqnarray}
which is a positive semidefinite biquadratic form that is not an SOS of bilinear forms \cite{Choi-75}.

Suppose $\lambda\in \mathbb{R}$ is an $M$-eigenvalue of $\mathscr{A}$ and $x,y\in \mathbb{R}^3\setminus \{0\}$ are the eigenvectors of $\mathscr{A}$ associated with the $M$-eigenvalue $\lambda$, then it follows from the definition of the $M$-eigenvalue that
\begin{eqnarray}\label{e-exam1-1}
\left\{\begin{array}{l}
x_1y_1^2-(x_2y_1y_2+x_3y_3y_1)+2x_1y_2^2=\lambda x_1,\\
x_2y_2^2-(x_1y_1y_2+x_3y_2y_3)+2x_2y_3^2=\lambda x_2,\\
x_3y_3^2-(x_2y_2y_3+x_1y_3y_1)+2x_3y_1^2=\lambda x_3,\\
x_1^2y_1-(x_1x_2y_2+x_3x_1y_3)+2x_3^2y_1=\lambda y_1,\\
x_2^2y_2-(x_1x_2y_1+x_2x_3y_3)+2x_1^2y_2=\lambda y_2,\\
x_3^2y_3-(x_2x_3y_2+x_3x_1y_1)+2x_2^2y_3=\lambda y_3,\\
x_1^2+x_2^2+x_3^2=1,\;\; y_1^2+y_2^2+y_3^2=1.
\end{array}\right.
\end{eqnarray}
From (\ref{e-exam1-1}) it is not difficult to show that $\lambda=0$ is an $M$-eigenvalue of $\mathscr{A}$, and
$$
((x_1,x_2,x_3),(y_1,y_2,y_3))=((1,0,0),(0,0,1)),\;\; ((0,1,0),(1,0,0)),\;\; ((0,0,1),(0,1,0))
$$
are the eigenvectors of $\mathscr{A}$ associated with the $M$-eigenvalue $\lambda=0$. Furthermore, $\lambda=0$ is the smallest $M$-eigenvalue of $\mathscr{A}$ since the biquadratic form $f(x,y)$ in (\ref{e-exam1-p-2}) is positive semidefinite.

The numerical result is listed in Table 1.
\begin{table}[ht] \label{table1}
\caption{The numerical results of the problem in Example \ref{exam1}}
\begin{center}
\tabcolsep 16.80pt
\begin{tabular}[c]{| c | c | c | c | c | c |}

\hline
  $r$& $s$& ITER&  CPU(s)& OPT              & VOL  \\ \hline

  0&   0&   9&     1.313& -1.5243441995e-09&  2.7246124908e-12 \\ \hline
  1&   0&  11&     1.719& -2.6041279810e-10&  5.7978715466e-12 \\ \hline
  0&   1&  11&     1.813& -2.6041187999e-10&  5.7989670655e-12\\ \hline
  1&   1&  13&     4.594& -2.4538473770e-10&  1.5420658730e-11 \\ \hline
  2&   1&  14&     8.797& -2.3617078649e-10&  1.3527256926e-11 \\ \hline
  1&   2&  14&     8.859& -2.3616477957e-10&  1.3665548633e-11 \\ \hline
  2&   2&  16&    29.344& -2.8204474346e-10&  1.0935180661e-11 \\ \hline
\end{tabular}
\end{center}
\end{table}

In Table 1, ``ITER" means the iteration number of the SDP solver, ``CPU(s)" means the total time in seconds spent for both setting up the problem and solving it, ``OPT" means the approximation value computed, and ``VOL" means the norm of the violation of the constraints of the approximation solution. From this table, we see that the method can find a good approximation solution even with the zero order relaxation.

\begin{Example}\label{exam2}
Let $\mathscr{A}=(a_{ijkl})\in \mathbb{T}_{4,3}$ be given by
\begin{eqnarray}\label{e-exam2-t-1}
\begin{array}{l}
a_{1111}=a_{1122}=a_{1133}=a_{2211}=a_{2222}=a_{2233}=a_{3311}=1,\\
a_{3322}=3,\; a_{3333}=3,\; a_{3323}=a_{3332}=-1,\; \mbox{\rm other}\; a_{ijkl}=0.
\end{array}
\end{eqnarray}
\end{Example}

The biquadratic form defined by this tensor is
\begin{eqnarray}\label{e-exam2-p-2}
%\begin{array}{rcl}
f(x,y)=x_1^2y_1^2 + x_1^2y_2^2 + x_1^2y_3^2 +x_2^2y_1^2 + x_2^2y_2^2 + x_2^2y_3^2+ x_3^2y_1^2+3x_3^2y_2^2 - 2x_3^2y_2y_3 + 3x_3^2y_3^2.
%\end{array}
\end{eqnarray}

Suppose $\lambda\in \mathbb{R}$ is an $M$-eigenvalue of $\mathscr{A}$ and $x,y\in \mathbb{R}^3\setminus \{0\}$ are the eigenvectors of $\mathscr{A}$ associated with the $M$-eigenvalue $\lambda$, then it follows from the definition of the $M$-eigenvalue that
\begin{eqnarray}\label{e-exam2-1}
\left\{\begin{array}{lll}
x_1y_1^2+x_1y_2^2+x_1y_3^2&=&\lambda x_1,\\
x_2y_1^2+x_2y_2^2+x_2y_3^2&=&\lambda x_2,\\
x_3y_1^2+3x_3y_2^2-2x_3y_2y_3+3x_3y_3^2&=&\lambda x_3,\\
x_1^2y_1+x_2^2y_1+x_3^2y_1&=&\lambda y_1,\\
x_1^2y_2+x_2^2y_2+3x_3^2y_2-x_3^2y_3&=&\lambda y_2,\\
x_1^2y_3+x_2^2y_3-x_3^2y_2+3x_3^2y_3&=&\lambda y_3,\\
x_1^2+x_2^2+x_3^2&=&1,\\
y_1^2+y_2^2+y_3^2&=&1.
\end{array}\right.
\end{eqnarray}
From the first three equalities and last two equalities in (\ref{e-exam2-1}), we can obtain that
$$
\lambda=1+x_3^2y_2^2+(x_3y_2-x_3y_2)^2+x_3^2y_3^2\geq 1,
$$
which implies that all $M$-eigenvalues of $\mathscr{A}$ given in (\ref{e-exam3-t-1}) are greater than or equal to $1$. Moreover, it is easy to see that $(\lambda_*, x^*, y^*)=(1, (1,0,0)^\top, (1,0,0)^\top)$ is a solution to (\ref{e-exam2-1}). Therefore, $\lambda_*=1$ is the smallest $M$-eigenvalues of $\mathscr{A}$ given in (\ref{e-exam3-t-1}).

The numerical result is listed in Table 2.
\begin{table}[ht] \label{table1}
\caption{The numerical results of the problem in Example \ref{exam2}}
\begin{center}
\tabcolsep 16.80pt
\begin{tabular}[c]{| c | c | c | c | c | c |}

\hline
  $r$& $s$& ITER&  CPU(s)& OPT              & VOL  \\ \hline
  0&   0&  10&     1.219&  9.9999999970e-01&  1.4472115766e-14 \\ \hline
  1&   0&  12&     1.734&  9.9999999963e-01&  1.8617918560e-11 \\ \hline
  0&   1&  11&     1.844&  9.9999999504e-01&  2.6733847183e-10 \\ \hline
  1&   1&  13&     4.906&  9.9999999954e-01&  1.7099952719e-11 \\ \hline
  2&   1&  14&     9.000&  9.9999999949e-01&  1.4822319679e-11 \\ \hline
  1&   2&  14&     8.938&  9.9999999812e-01&  3.8117223117e-11 \\ \hline
  2&   2&  16&    29.734&  9.9999999749e-01&  6.7937294113e-12 \\ \hline
\end{tabular}
\end{center}
\end{table}

In Table 2, symbols ``ITER", ``CPU(s)", ``OPT", ``VOL" are the same as in Example \ref{exam1}. From this table, we see that the method can find a good approximation solution even with the zero order relaxation.

\noindent {\bf Part 2}. In this part, by using some results obtained in Section 5, we compute the smallest eigenvalue of an unfolded matrix of some semi-paired symmetric tensor to judge whether the concerned paired symmetric tensor is positive definite or not.

\begin{Example}\label{exam3}
Let $\mathscr{A}=(a_{ijkl})\in \mathbb{T}_{4,3}$ be given by
\begin{eqnarray}\label{e-exam3-t-1}
\begin{array}{l}
a_{1111}=2,\; a_{2222}=2,\; a_{3333}=2,\; a_{1313}=a_{3113}=a_{1331}=a_{3131}=-\frac 12,\\
a_{2323}=a_{3223}=a_{2332}=a_{3232}=-\frac 12,\; \mbox{\rm other}\; a_{ijkl}=0.
\end{array}
\end{eqnarray}
\end{Example}

The biquadratic form defined by this tensor is
\begin{eqnarray}\label{e-exam3-p-2}
f(x,y)=2x_1^2y_1^2+2x_2^2y_2^2+2x_3^2y_3^2-2x_1x_3y_1y_3-2x_2x_3y_2y_3.
\end{eqnarray}
It is not difficult to obtain that $f(x,y)=x_1^2y_1^2+x_2^2y_2^2+(x_1y_1-x_3y_3)^2+ (x_2y_2-x_3y_3)^2$ and $f(x^*,y^*)=0$ where $x^*=(1,0,0)^\top$ and $y^*=(0,1,0)$. These demonstrate that the biquadratic form in (\ref{e-exam3-p-2}) is an SOS of bilinear forms but not positive definite.

Let $\mathscr{A}$ be given by (\ref{e-exam3-t-1}) and $\mathscr{B}$ defined by Remark \ref{remark-def-semi-p-s} be a semi-paired symmetric tensor, then
\begin{eqnarray*}
\begin{array}{l}
b_{1111}=2,\; b_{2222}=2,\; b_{3333}=2,\;  b_{1313}=d_{3131}=-1,\;
b_{2323}=b_{3232}=-1,\;  \mbox{\rm other}\; b_{ijkl}=0.
\end{array}
\end{eqnarray*}
Suppose $M\in \mathbb{R}^{9\times 9}$ defined by (\ref{e-matrix-A-0}) is the unfolded matrix of $\mathscr{B}$, then
\begin{eqnarray*}
m_{11}=2,\; m_{55}=2,\; m_{99}=2,\; m_{19}=m_{91}=-1,\; m_{59}=m_{95}=-1,\;  \mbox{\rm other}\; m_{ij}=0.
\end{eqnarray*}
It is easy to obtain that eigenvalues of $M$ are $0, 0, 0, 0, 0, 0, 0.5858, 2, 3.4142$, respectively, which implies that the matrix $M$ is positive semidefinite but not positive definite. Thus, by Corollary \ref{coro-sos-4} we obtain that the quadratic form (\ref{e-exam3-p-2}) is an SOS of bilinear forms but not positive definite.

\begin{Example}\label{exam4}
Let $\mathscr{A}=(a_{ijklpq})\in \mathbb{T}_{6,3}$ be given by
\begin{eqnarray}\label{e-exam4-t-1}
\begin{array}{l}
a_{111111}=a_{112233}=a_{221111}= a_{222222}=a_{331122}=a_{333333}=1,\\
a_{111213}=a_{111231}=a_{112113}=a_{112131}=-\frac 12,\;\; a_{221212}=a_{221221}=a_{222112}=a_{222121}=-\frac 12,\\
 a_{331323}=a_{331332}=a_{333123}=a_{333132}=-\frac 12,\;\; \mbox{\rm other}\; a_{ijklpq}=0.
\end{array}
\end{eqnarray}
\end{Example}

The homogeneous polynomial defined by this tensor is
\begin{eqnarray}\label{e-exam4-p-2}
\begin{array}{rcl}
f(x,y,z)&=&x_1^2y_1^2z_1^2+x_1^2y_2^2z_3^2+x_2^2y_1^2z_1^2+x_2^2y_2^2z_2^2+x_3^2y_1^2z_2^2+x_3^2y_3^2z_3^2\\
& &-2(x_1^2y_1y_2z_1z_3+x_2^2y_1y_2z_1z_2+x_3^2y_1y_3z_2z_3),
\end{array}
\end{eqnarray}
It is not difficult to obtain that $f(x,y,z)=(x_1y_1z_1-x_1y_2z_3)^2+ (x_2y_1z_1-x_2y_2z_2)^2+(x_3y_1z_2-x_3y_3z_3)^2$ and $f(x^*,y^*,z^*)=0$ where $x^*=(0,0,1)^\top, y^*=(0,1,0)$ and $y^*=(1,0,0)$. These demonstrate that the homogeneous polynomial in (\ref{e-exam4-p-2}) is an SOS of bilinear forms but not positive definite.

Let $\mathscr{A}$ be given by (\ref{e-exam4-t-1}) and $\mathscr{B}$ defined by Remark \ref{remark-def-semi-p-s-6} be a semi-paired symmetric tensor, then
\begin{eqnarray*}
\begin{array}{l}
b_{111111}=b_{112233}=b_{221111}= b_{222222}=b_{331122}=b_{333333}=1,\\
b_{111213}=b_{112131}=-1,\;\; b_{221212}=b_{222121}=-1,\; b_{331323}=b_{333132}=-1,\;\; \mbox{\rm other}\; b_{ijklpq}=0.
\end{array}
\end{eqnarray*}
Suppose $N\in \mathbb{R}^{27\times 27}$ defined by (\ref{e-matrix-N-0}) is the unfolded matrix of $\mathscr{B}$, then
\begin{eqnarray*}
\begin{array}{l}
n_{11}=n_{66}=n_{10,10}=n_{14,14}=n_{20,20}=n_{27,27}=1,\\
n_{16}=n_{61}=n_{10,14}=n_{14,10}=n_{20,27}=n_{27,20}=-1,\;  \mbox{\rm other}\; n_{ij}=0.
\end{array}
\end{eqnarray*}
It is easy to obtain that eigenvalues of $N$ are $\underbrace{0,0,\cdots, 0}_{24}, 2, 2, 2$, respectively, which implies that the matrix $N$ is positive semidefinite but not positive definite. Thus, by Corollary \ref{coro-sos-6} we obtain that the homogeneous polynomial in (\ref{e-exam4-p-2}) is an SOS of trilinear forms but not positive definite.

\section{Concluding Remarks}

In this paper, we investigated positive definiteness of the fourth order three dimensional and sixth order three dimensional paired symmetric (elasticity) tensor. With the help of the $M$-eigenvalue we gave several necessary and/or sufficient conditions under which a fourth order three dimensional or sixth order three dimensional paired symmetric (elasticity) tensor is positive definite. By introducing semi-paired symmetric tensor and the unfolded matrix of that tensor, we showed that positive definiteness of a fourth order three dimensional and sixth order three dimensional paired symmetric (elasticity) tensor can be obtained by investigating positive definiteness of some unfolded matrix of a semi-paired symmetric tensor of that tensor. We also proposed a necessary and sufficient condition under which a fourth order three dimensional or sixth order three dimensional paired symmetric (elasticity) tensor is positive semidefinite. These results can be extended to the case of higher order higher dimensional (strongly) paired symmetric tensors. We reported preliminary numerical results which confirm our theoretical findings. 

We derived several necessary and sufficient conditions for positive definiteness of a sixth order three dimensional paired symmetric tensor in Section 4, where the most main condition needs to judge positive definiteness of a fourth order paired symmetric tensor, an eighth order bi-block symmetric tensor and a twelfth order bi-block symmetric tensor. It is well known that it is difficult to solve the polynomial optimization problem with the involved polynomials being defined by higher order higher dimensional tensors. Thus, one of further issues is how to design effective methods to judge positive definiteness of higher order (strongly) paired symmetric tensors or higher order elasticity tensors. Moreover, in this paper, we mainly gave some analysis for positive definiteness of the fourth order three dimensional and sixth order three dimensional paired symmetric (elasticity) tensor. It is worth investigating effective algorithms to check positive definiteness of paired symmetric (elasticity) tensors arised from practical mechanical problems.

\end{document}